\documentclass{article}

\usepackage{indentfirst}
\usepackage{PRIMEarxiv}
\usepackage{amsmath,amssymb,amsthm}
\usepackage[ruled,vlined]{algorithm2e}
\usepackage[utf8]{inputenc} 
\usepackage[T1]{fontenc}    
\usepackage{hyperref}       
\usepackage{url}            
\usepackage{booktabs}       
\usepackage{amsfonts}       
\usepackage{nicefrac}       
\usepackage{microtype}      
\usepackage{lipsum}
\usepackage{graphicx}
\usepackage{color}
\usepackage{subfigure}
\usepackage{authblk}
\graphicspath{{media/}}     

\title{Symplectic Discretization Approach for Developing New Proximal Point Algorithm
\thanks{This work was partially supported by grant 12288201 from NSF of China.
}}

\author[1]{Ya-xiang Yuan}
\author[ 1,2]{Yi Zhang\thanks{Corresponding author: zhangyi2020@lsec.cc.ac.cn}  }
\affil[1]{Institute of Computational Mathematics and Scientific/Engineering Computing, Academy of Mathematics and Systems
	Science, Chinese Academy of Sciences, Beijing 100190, China}
\affil[2]{University of Chinese Academy of Sciences, Beijing 100049, China}

\newcommand{\argmin}{\mathop{\arg\min}}
\newcommand{\inner}[1]{\left\langle#1\right\rangle}
\newcommand{\norm}[1]{\left\|#1\right\|}
\newcommand{\X}{\overset{\cdot}{X}}
\newcommand{\XX}{\overset{\cdot\cdot}{X}}
\newcommand{\E}{\mathcal{E}}
\newcommand{\EE}{\overset{\cdot}{\mathcal{E}}}
\newcommand{\Z}{\overset{\cdot}{Z}}

\newcommand{\G}{\mathcal{G}}
\newcommand{\hh}{\mathcal{H}}
\newcommand{\dist}{\mathrm{dist}}
\newcommand{\zero}{\mathrm{zero}}
\newcommand{\dd}{\mathrm{d}}

\newtheorem{theorem}{Theorem}

\newtheorem{lemma}{Lemma}
\newtheorem{proposition}{Proposition}
\begin{document}
\maketitle

\begin{abstract}
The rapid advancements in high-dimensional statistics and machine learning have increased the use of first-order methods. Many of these methods can be regarded as instances of the proximal point algorithm. Given the importance of the proximal point algorithm, there has been growing interest in developing its accelerated variants. However, some existing accelerated proximal point algorithms exhibit oscillatory behavior, which can impede their numerical convergence rate. In this paper, we first introduce an ODE system and demonstrate its \( o(1/t^2) \) convergence rate and weak convergence property. Next, we apply the Symplectic Euler Method to discretize the ODE and obtain a new accelerated proximal point algorithm, which we call the Symplectic Proximal Point Algorithm. The reason for using the Symplectic Euler Method is its ability to preserve the geometric structure of the ODEs. Theoretically, we demonstrate that the Symplectic Proximal Point Algorithm achieves an \( o(1/k^2) \) convergence rate and that the sequences generated by our method converge weakly to the solution set. Practically, our numerical experiments illustrate that the Symplectic Proximal Point Algorithm significantly reduces oscillatory behavior, leading to improved long-time behavior and faster numerical convergence rate.
\end{abstract}


\section{Introduction}
With the rapid advancements of high-dimensional statistics and machine learning, first-order methods have become widely utilized for training models. Many of these first-order methods have been found to be closely related to a fundamental root-finding algorithm, known as the proximal point algorithm (PPA). For example, the augmented Lagrangian method (ALM), the alternating direction method of multipliers (ADMM), and the Chambolle-Pock method (also known as the primal-dual hybrid gradient method, PDHG) have been proven to be instances of PPA, as discussed in \cite{rockafellar78, eckstein92, chambolle11, chambolle16, connor20}. In addition, several first-order methods for solving variational inequality problem, such as the extra-gradient method and the optimistic gradient descent ascent method, have also been found to be closely related to PPA, as detailed in \cite{mokhtari20}.

The convergence rate of PPA has been proven to be \( O(1/k) \) \cite{baillon96}, which can lead to slow numerical performance in practical applications. Therefore, accelerating the PPA is crucial for improving its practical efficiency. Nesterov's accelerated gradient method (NAG), which demonstrates an \( O(1/k^2) \) convergence rate and has the same computational complexity as vanilla gradient descent, has inspired research into accelerated variants of the PPA. Recently, an accelerated technique known as Halpern's iteration \cite{halpern67} has gained increasing attention. \cite{lieder21} and \cite{qi21} have shown that the convergence rate of Halpern's iteration is \( O(1/k^2) \). Additionally, \cite{kim21} introduced a general analytical framework for analyzing the PPA and demonstrated the \( O(1/k^2) \) convergence rate of Halpern's iteration. The work in \cite{tran24} also studied the connection between Halpern's iteration and the NAG method. Due to the fascinating theoretical results of Halpern's iteration, the study of its applications has grown. The works \cite{zhang24, liang24} used Halpern's iteration to develop new algorithms for solving specific optimization problems. Furthermore, since \cite{mokhtari20} suggested that the extra-gradient method and the optimistic gradient descent ascent method are closely related to the PPA, Halpern's iteration has been employed to accelerate these methods, as shown in \cite{cai24, lee21, park22, tran21, yoon21}.

With the growing research on the connections between optimization algorithms and ordinary differential equations (ODEs), deriving new accelerated PPA frameworks has become feasible. Initially, \cite{su16} linked the NAG method with a second-order ODE and introduced the Lyapunov function technique to demonstrate the $O(1/k^2)$ convergence rate for the NAG method. Recent studies have increasingly utilized Lyapunov functions to establish the convergence properties of optimization algorithms, as exemplified in \cite{attouch16,chen23,li2022proximalsubgradientnormminimization,shi19}. Given the theoretical depth of ODEs, several studies have focused on deriving accelerated versions of the PPA from the discretization of ODEs. In recent work \cite{attouch20,attouch21,attouch192,bot23}, researchers first proposed several second-order ODEs and then studied convergence properties of these ODEs. By applying the explicit or implicit difference schemes to these continuous ODEs for discretization, they successfully developed corresponding accelerated PPA variants.

\subsection{Existing Results}
\label{sec:1.1}
Let $\hh$ be a real Hilbert space equipped with the inner product $\inner{\cdot, \cdot}$ and the corresponding norm $\norm{\cdot}$ and let $A: \hh\to 2^{\hh}$ be a set-valued operator operator. $A$ is monotone if the following inequality holds:
\[
\inner{u-v, x-y}\geqslant 0,\quad\forall x, y\in\hh,  u\in A(x), v\in A(y).
\]
The monotone operator $A$ is maximal if the graph of $A$, defined as $\{(x,u)\in \hh\times\hh| u\in A(x)\}$, is not properly included by the graph of other monotone operator.

In this paper, we consider the following monotone inclusion problem:
\begin{equation}
	\label{eq:main}
	0\in A(x),
\end{equation}
where $A$ is maximally monotone. Let $\zero(A):=\{x\in\hh|0\in A(x)\}$ be the solution set of \eqref{eq:main}, assuming that $\zero(A)$ is non-empty. The PPA for solving \eqref{eq:main} is given as follows:\begin{equation}
	\label{eq:ppa}
	x_{k+1}=J_{cA}(x_k),\quad c>0,
\end{equation}
which can be interpreted as applying Banach-Picard iteration to the proximal point operator $J_{cA}:=(I+cA)^{-1}$. The parameter $c$ refers to the index of the proximal point operator $J_{cA}$. Here we assume that $J_{cA}$ can be easily calculated for some specific $c$. For simplicity, we assume $c=1$ for the rest of this paper. The convergence rate of \eqref{eq:ppa} is
\[
\norm{x_{k+1}-x_k}^2\leqslant\frac{1}{1+k}\cdot\dist(x_0,\zero(A))^2,
\]
where $\dist(x_0,\zero(A))=\inf\limits_{x\in\zero(A)}\norm{x_0-x}$. We refer to the lecture note \cite{baillon96} and the textbooks \cite{bauschke17} and \cite{ryu20} for the proof of above result. The term $\norm{x_{k+1}-x_k}^2$ can be used to characterize convergence rate of PPA is due to the property $x_k-x_{k+1}\in A(x_{k+1})$. Recently, \cite{gu20} showed a tighter rate of PPA:\[
\norm{x_{k+1}-x_k}^2\leqslant\cfrac{1}{(k+1)\left(1+\cfrac 1k\right)^k}\cdot\dist(x_0,\zero(A))^2.
\]

The convergence of the sequence generated by the PPA is also of concern. In this paper, the notation \( x_k \rightharpoonup x \) denotes that the sequence \(\{x_k\}\) converges weakly to \( x \), and the notation \( x_k \to x \) represents that the sequence \(\{x_k\}\) converges strongly to \( x \). The weak convergence property of the PPA relies on the following propositions:

\begin{proposition}
	\label{prop:discreteopial}
	Let \(\{x_k\}\) be a sequence in \(\mathcal{H}\) and let \( S \) be a nonempty subset of \(\mathcal{H}\). Suppose that, for every \( x \in S \), \(\norm{x_k - x}\) converges and that every weak sequential cluster point of \(\{x_k\}\) belongs to \( S \). Then \(\{x_k\}\) converges weakly to a point in \( S \).
\end{proposition}

Proposition \ref{prop:discreteopial} can be seen as a discrete version of Opial's Lemma \cite{opial67}.

\begin{proposition}
	\label{prop:continuousopial}
	Let \( X(t): [t_0, +\infty) \to \mathcal{H} \) and let \( S \) be a nonempty subset of \(\mathcal{H}\). Suppose that, for every \( x \in S \), \(\norm{X(t) - x}\) converges and that every weak sequential cluster point of \( X(t) \) belongs to \( S \). Then \( X(t) \) converges weakly to a point in \( S \).
\end{proposition}

\begin{proposition}
	\label{prop:nonexpansiveweak}
	Let \( A \) be a maximally monotone operator. If there exist two sequences \(\{x_k\}\) and \(\{u_k\}\) such that \( u_k \in A(x_k) \), \( x_k \rightharpoonup x \), and \( u_k \to u \), then \( u \in A(x) \).
\end{proposition}

Proposition \ref{prop:discreteopial} and Proposition \ref{prop:nonexpansiveweak} correspond to Lemma 2.47 and Proposition 20.38 (ii) in \cite{bauschke17}. One can prove the weak convergence property of the PPA by using Proposition \ref{prop:discreteopial} and Proposition \ref{prop:nonexpansiveweak}, along with the facts that \( x_k - x_{k+1} \in A(x_{k+1}) \) and \( x_k - x_{k+1} \to 0 \). For further theory and applications of the PPA, we refer to \cite{bauschke17, combettes18, parikh14, ryu20}.

Let $T:\hh\to\hh$ be a $1-$Lipschitz operator. The Halpern's iteration proposed in \cite{halpern67} can be described as follows:\begin{equation}
	\label{eq:halperniteration}
	x_{k+1}=\alpha_k x_0+(1-\alpha_k)T(x_k),\quad \alpha_k\in(0,1).
\end{equation}
By setting \(\alpha_k = \dfrac{1}{k+2}\), \(T = 2J_A - I\), and omitting \(x_0\) from the iteration rule, we can derive the following accelerated proximal point algorithm from \eqref{eq:halperniteration}:
\begin{equation}
	\begin{split}
		y_{k+1} & = J_{A}(x_k), \\
		x_{k+1} & = y_{k+1}+\frac{k}{k+2}(y_{k+1}-y_k)-\frac{k}{k+2}(y_k-x_{k-1}),
	\end{split}
	\label{eq:acceleratedproximal}
\end{equation}
where $y_0=x_0$. In \cite{kim21} and \cite{ryu20}, the convergence rate of \eqref{eq:acceleratedproximal} has been proven to be\[
\norm{y_{k+1}-x_k}^2\leqslant\frac{1}{(k+1)^2}\cdot\dist(x_0,\zero(A))^2.
\]
Also by Theorem 30.1 in \cite{bauschke17}, the sequence $\{x_k\}$ generated by \eqref{eq:acceleratedproximal} strongly converges to the projection of $x_0$ onto $\zero(A)$.

Nowadays, there are several accelerated variants of the PPA derived from ODEs. Initially, \cite{attouch192, attouch20, attouch21} used the Yosida approximation of $A$ to propose some second-order ODEs, and then obtained accelerated PPAs by properly discretizing these ODEs. However, the index of the proximal point operator in the resulting accelerated PPAs changes as iterations progress, which may limit their applications in deriving some accelerated convex optimization algorithms (see \cite{eckstein92} for discussion). More recently, \cite{bot23} proposed a new accelerated PPA framework. They were inspired by \cite{bot232} and studied the following ODE:\[
\overset{\cdot\cdot}{X}+\frac{\alpha}{t}\X+\frac{\mathrm{d}}{\mathrm{d}t}(I-J_A)+\frac{\alpha}{2t}(I-J_A)=0,
\]
where $\alpha>2$ is a constant. Then, they used the explicit difference scheme to discretize the above ODE and obtained the following Fast K-M algorithm:\begin{equation}
	\label{eq:fastkm}
	\begin{split}
		x_{k+1}=&\ \left(1-\frac{s\alpha}{2(k+\alpha)}\right)x_k+\frac{(1-s)k}{k+\alpha}(x_k-x_{k-1})\\
		&\ +\frac{s\alpha}{2(k+\alpha)}J_{A}(x_k)+\frac{sk}{k+\alpha}(J_{A}(x_k)-J_{A}(x_{k-1})),
	\end{split}
\end{equation}
where $s>0$ is a constant. They proved that the convergence rate of \eqref{eq:fastkm} is $o(1/k^2)$ and the sequence $\{x_k\}$ generated by \eqref{eq:fastkm} converges weakly to a point in $\zero(A)$.
\subsection{Our Motivations and Contributions}
First, we introduce the motivation for our study. In \cite{bot23}, numerical experiments showed that \eqref{eq:fastkm} suffers from oscillation phenomenon, which may slow down the numerical convergence rate of \eqref{eq:fastkm}. Based on this observation, we develop the \emph{Symplectic Proximal Point Algorithm (SPPA)}. The details of this development are as follows.

First, we need to construct an ODE system and analyze its convergence properties, which will provide insight into the desired algorithm. To obtain an ODE that exhibits an \(o(1/t^2)\) convergence rate, we draw inspiration from the ODEs in \cite{bot23} and \cite{yi23}. Based on these insights, we propose the following ODE:
\begin{equation}
	\label{eq:acceleratedode}
	\begin{split}
		Z &= \frac{t}{r}\dot{X} + \left(1 + \frac{t}{r}\right)A(X) + X,\\
		\dot{Z} &= -\frac{C}{r}A(X),\\
		Z(0) &= X(0) = x_0,
	\end{split}
\end{equation}
where $r$ and $C$ are two positive constants. Here, \(A\) is assumed to be single-valued and \(A(X)\) is differentiable with respect to \(t\). We show that the convergence rate of the solution trajectory to \eqref{eq:acceleratedode} is \(o(1/t^2)\), and that the solution trajectory converges weakly to $\zero(A)$.

Next, we need to discretize \eqref{eq:acceleratedode}. The discretization method we used is the \emph{Symplectic Euler Method}, which is why we refer to our method as the Symplectic Proximal Point Algorithm. Our motivation for using the Symplectic Euler Method comes from the work in \cite{shi19,shi22}, where the authors showed that the NAG method can be viewed as applying the symplectic method to high-resolution ODEs formulated in phase space representation. Since the Symplectic Euler Method can preserve geometric structure of the ODEs and has a simple iteration rule, we use the Symplectic Euler Method to derive our algorithm. 

First, we briefly introduce the Symplectic Euler Method. Consider the following Hamiltonian system:
\begin{equation}
	\begin{split}
		\dot{p} &= -\frac{\partial H(p,q)}{\partial q},\\
		\dot{q} &= \frac{\partial H(p,q)}{\partial p}.
	\end{split}
\end{equation}
The Symplectic Euler Method is given as follows:
\begin{equation}
	\label{eq:symplecticeuler}
	\begin{split}
		p_{k+1} &= p_k - s\frac{\partial H(p_k, q_{k+1})}{\partial q},\\
		q_{k+1} &= q_k + s\frac{\partial H(p_k, q_{k+1})}{\partial p},
	\end{split}
	\quad \text{or} \quad
	\begin{split}
		p_{k+1} &= p_k - s\frac{\partial H(p_{k+1}, q_k)}{\partial q},\\
		q_{k+1} &= q_k + s\frac{\partial H(p_{k+1}, q_k)}{\partial p}.
	\end{split}
\end{equation}
For more information about symplectic methods, we refer to the textbooks \cite{hairer06,feng10}. An introduction to the Symplectic Euler Method can be found in Section I.1.2 of \cite{hairer06}.

\begin{figure}[ht]
	\centering
	\subfigure{\includegraphics[width=.3\textwidth]{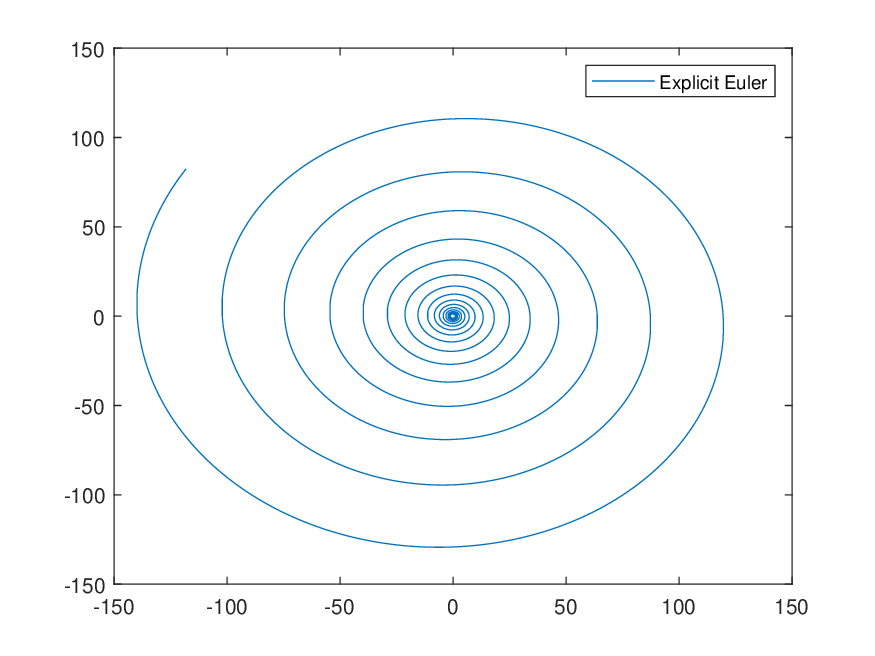}}
	\subfigure{\includegraphics[width=.3\textwidth]{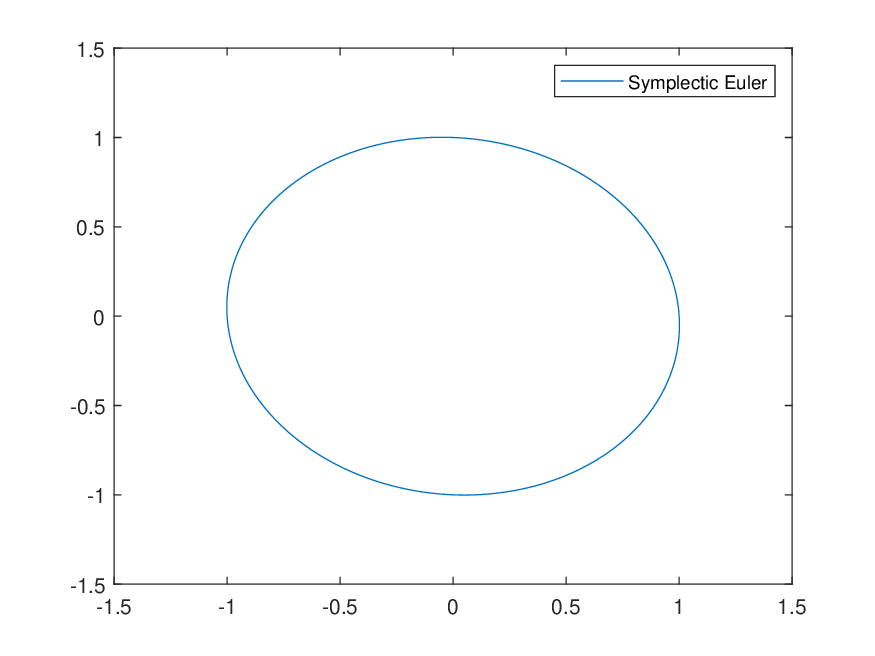}}
	\subfigure{\includegraphics[width=.3\textwidth]{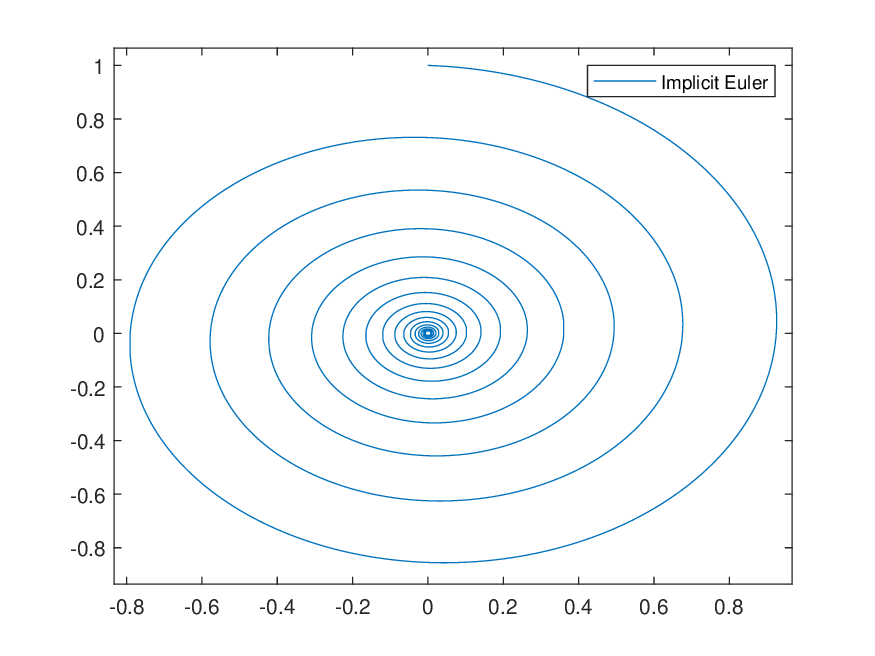}}
	\caption{An example illustrates that the Symplectic Euler Method can preserve geometric structure of the ODEs. Here, we apply the Explicit Euler Method, Symplectic Euler Method, and Implicit Euler Method to the Hamiltonian system \(\dot{p} = q\), \(\dot{q} = -p\) with initial condition \(p(0) = 0\) and \(q(0) = 1\). The exact solution to this Hamiltonian system is \(p = \sin t\) and \(q = \cos t\). As shown in the above figures, the Symplectic Euler method preserves geometric structure of the solution.}
\end{figure}

By applying (\ref{eq:symplecticeuler}) to discretize (\ref{eq:acceleratedode}), we obtain the following iteration rule:
\begin{equation}
	\begin{split}
		z_0&=x_0,\\
		z_k &= \frac{k}{r}(x_{k+1} - x_k) + \left(1 + \frac{k}{r}\right)A(x_{k+1}) + x_{k+1},\\
		z_{k+1} &= z_k - \frac{C}{r}A(x_{k+1}).
	\end{split}
	\label{eq:SPPAtemp}
\end{equation}
However, (\ref{eq:SPPAtemp}) requires \(A\) to be single-valued. Therefore, we reformulate the algorithm using the definition of proximal point operator. Specifically, the update rule of \(x_{k+1}\) can be transformed into:
\[
\begin{split}
	\tilde{x}_{k+1} &= \frac{k}{k+r}x_k + \frac{r}{k+r}z_k,\\
	x_{k+1} &= J_A(\tilde{x}_{k+1}).
\end{split}
\]
Moreover, since \(\tilde{x}_{k+1} - x_{k+1} \in A(x_{k+1})\), the update rule for \(z_{k+1}\) can be transformed into:
\[
z_{k+1} = z_k + \frac{C}{r}(x_{k+1} - \tilde{x}_{k+1}).
\]
The resulting algorithm is detailed in Algorithm \ref{al:sppa}.

\begin{algorithm}[H]
	\caption{Symplectic Proximal Point Algorithm, SPPA}
	\label{al:sppa}
	Initialize $x_0, z_0=x_0$\;
	\For{$k=0, 1, \cdots$}{
		$\tilde{x}_{k+1}=\dfrac{k}{k+r}x_k+\dfrac{r}{k+r}z_k$\;
		$x_{k+1}=J_A(\tilde{x}_{k+1})$\;
		$z_{k+1}=z_k+\dfrac{C}{r}(x_{k+1}-\tilde{x}_{k+1})$\;	
	}
\end{algorithm}

We also demonstrate that Algorithm \ref{al:sppa} exhibits an \(o(1/k^2)\) convergence rate, and the sequences generated by Algorithm \ref{al:sppa} converge weakly to the solution set without assuming that \(A\) is single-valued. Furthermore, our numerical experiments show that the oscillation phenomenon is less significant compared to some existing accelerated PPAs, and the numerical performance of SPPA is better than that of some existing accelerated PPAs.

\section{Continuous Time Perspective}
\label{sec:2}

In this section, we investigate the convergence rate and weak convergence property of \eqref{eq:acceleratedode}. To prove an \(o(1/t^2)\) convergence rate for \eqref{eq:acceleratedode}, we need to establish an \(O(1/t^2)\) convergence rate and integral rate for \eqref{eq:acceleratedode}. To achieve these results, we propose the following Lyapunov function, inspired by those given in \cite{yoon21} and \cite{yi23}:
\begin{equation}
	\label{eq:lyapunovcontinuous}
	\E(t) = \frac{t(t + r)}{2r^2}\norm{A(X(t))}^2 + \frac{t}{r}\inner{A(X(t)), X(t) - Z(t)} + \frac{r - 1}{2C}\norm{Z(t) - x^*}^2,
\end{equation}
where $x^*\in\zero(A)$. For simplicity, when proving the lemma and the theorem presented in this section, we use $A(X)$, $X$, and $Z$ to represent $A(X(t))$, $X(t)$, and $Z(t)$, respectively. The notations $A(X(t))$, $X(t)$ and $Z(t)$ are used to ensure the correct representation of our results. Next, we will examine its non-negativity and derivative.

\begin{lemma}
	\label{lem:lyapunovcontinousnegative}
	Let $A$ be a single-valued maximally monotone operator, and let $(X(t),Z(t))$ be the solution trajectory to \eqref{eq:acceleratedode}. The Lyapunov function \eqref{eq:lyapunovcontinuous} can be represented as\begin{align*}
		\E(t)= &\ \left[\left(1-\frac{C}{r-1}\right)\frac{t^2}{2r^2}+\frac{t}{2r}\right]\norm{A(X(t))}^2+\frac{t}{r}\inner{A(X(t)),X(t)-x^*}\\
		&\ +\frac{1}{2}\norm{\sqrt{\frac{Ct^2}{r^2(r-1)}}A(X(t))+\sqrt{\frac{r-1}{C}}(x^*-Z(t))}^2.
	\end{align*}
	In addition, if $0<C\leqslant r-1$, $\E(t)$ can be represented by a sum of three non-negative terms.
\end{lemma}
\begin{proof}
	 First, we have
		\[
		\frac{t}{r}\inner{A(X),X-Z}=\frac{t}{r}\inner{A(X),X-x^*}+\frac{t}{r}\inner{A(X),x^*-Z}.
		\]
	By the equality\begin{equation*}
		\inner{a,b}=\frac{1}{2}\left(\norm{a+b}^2-\norm{a}^2-\norm{b}^2\right),\quad\forall a, b\in\hh,
	\end{equation*}
	we have\begin{align*}
		\frac{t}{r}\inner{A(X),x^*-Z}	= &\ \frac{1}{2}\norm{\sqrt{\frac{Ct^2}{r^2(r-1)}}A(X)+\sqrt{\frac{r-1}{C}}(x^*-Z)}^2\\
		&\ -\frac{Ct^2}{2r^2(r-1)}\norm{A(X)}^2-\frac{r-1}{2C}\norm{Z-x^*}^2 
	\end{align*}
	Thus, we have\begin{align*}
		\E(t)= &\ \left[\left(1-\frac{C}{r-1}\right)\frac{t^2}{2r^2}+\frac{t}{2r}\right]\norm{A(X)}^2+\frac{t}{r}\inner{A(X),X-x^*}\\
		&\ +\frac{1}{2}\norm{\sqrt{\frac{Ct^2}{r^2(r-1)}}A(X)+\sqrt{\frac{r-1}{C}}(x^*-Z)}^2.
	\end{align*}
	Also, if $C\leqslant r-1$, $1-\dfrac{C}{r-1}\geqslant 0$. Thus $\left(1-\dfrac{C}{r-1}\right)\dfrac{t^2}{2r^2}+\dfrac{t}{2r}\geqslant 0$. In conclusion, $\E(t)$ can be represented by a sum of three non-negative terms.
\end{proof}
\begin{lemma}
	\label{lem:lyapunovcontinuous}
	Let $A$ be a single-valued maximally monotone operator, and let $(X(t), Z(t))$ be the solution trajectory to \eqref{eq:acceleratedode}. If $A(X)$ is differentiable respect to $t$, we have
	\begin{align*}
		\EE(t)=&\ -\frac{t^2}{r^2}\inner{\frac{\dd}{\dd t}A(X(X(t))),\X(t)}-\left(\frac{r-1-C}{r^2}t+1-\frac{1}{2r}\right)\norm{A(X(t))}^2\\
		&\ -\frac{r-1}{r}\inner{A(X(t)),X(t)-x^*}.
	\end{align*}
	
\end{lemma}
\begin{proof}
	Direct calculation shows that\begin{align*}
		 \EE(t)
		= &\ \frac{t(t+r)}{r^2}\inner{\frac{\dd}{\dd t}A(X),A(X)}+\frac{2t+r}{2r^2}\norm{A(X)}^2+\frac{1}{r}\inner{A(X),X-Z}\\
		&\ +\frac{t}{r}\inner{\frac{\dd}{\dd t}A(X),X-Z}+\frac{t}{r}\inner{A(X),\X-\Z}+\frac{r-1}{C}\inner{\Z,Z-x^*}\\
		=&\ -\frac{t^2}{r^2}\inner{\frac{\dd}{\dd t}A(X),\X}+\frac{2t+r}{2r^2}\norm{A(X)}^2+\frac{1}{r}\inner{A(X),X-Z}\\
		&\ +\inner{A(X), Z-X-\left(1+\frac{t}{r}-\frac{Ct}{r^2}\right)A(X)}-\frac{r-1}{r}\inner{A(X),Z-x^*}\\
		=&\ -\frac{t^2}{r^2}\inner{\frac{\dd}{\dd t}A(X),\X}-\left(\frac{r-1-C}{r^2}t+1-\frac{1}{2r}\right)\norm{A(X)}^2\\
		&\ -\frac{r-1}{r}\inner{A(X),X-x^*}.
	\end{align*}
	
\end{proof}

By using previous two results, we can obtain the following theorem.
\begin{theorem}
	\label{thm:convergenceratecontinuous}
	Let $A$ be a single-valued maximally monotone operator, and let $(X(t),Z(t))$ be the solution to \eqref{eq:acceleratedode}. If  $A(X)$ is differentiable respect to $t$,	$r>1$, and $0<C\leqslant r-1$,  then $\E(t)$ is non-increasing and non-negative. Also, the following inequalities hold: \begin{align}
		\norm{A(X(t))}^2 & \leqslant\frac{r^2(r-1)^2}{[C(r-1)-C^2]t^2+Cr(r-1)t}\cdot \dist(x_0,\zero(A))^2,\label{eq:t2bound} \\
		\inner{A(X(t)),X(t)-x^*} &\leqslant\frac{r^2-r}{2Ct}\norm{x_0-x^*}^2,\quad\forall x^*\in\zero(A).\label{eq:t1bound}
	\end{align}
	Additionally, we have
	\begin{align}
		\int_{0}^{\infty}\left(\frac{r-1}{r}t+\frac{1}{2r}\right)\norm{A(X(t))}^2\dd t&\leqslant\frac{r-1}{2C}\dist(x_0,\zero(A))^2,\label{eq:integralfinite}\\
		\int_{0}^{\infty}\left(1-\frac{1}{r}\right)\inner{A(X(t)),X(t)-x^*}\dd t&\leqslant\frac{r-1}{2C}\norm{x_0-x^*}^2,\quad\forall x^*\in\zero(A).\label{eq:integralfinite2} 
	\end{align}
\end{theorem}
\begin{proof}
	Because $A$ is monotone, we have that\[
	\inner{A(X(t+\Delta t))-A(X(t)),X(t+\Delta t)-X(t)}\geqslant 0,\quad\forall t\geqslant 0,\ \Delta t>0.
	\]
	By dividing $\Delta t^2$ on both sides of above inequality and letting $\Delta t\to 0$, we have\[
	\inner{\frac{\dd}{\dd t}A(X),\X}\geqslant 0.
	\]
	Owing to Lemma \ref{lem:lyapunovcontinuous} and above inequality, we have\[
	\EE(t)\leqslant-\left(\frac{r-1-C}{r^2}t+1-\frac{1}{2r}\right)\norm{A(X)}^2-\frac{r-1}{r}\inner{A(X),X-x^*}.
	\]
	Since $r>1$ and $C\leqslant r-1$, $\dfrac{r-1-C}{r^2}t+1-\dfrac{1}{2r}$ and $\dfrac{r-1}{r}$ are positive. Thus $\EE(t)\leqslant 0$, which implies non-increasing of $\E(t)$. Then by Lemma \ref{lem:lyapunovcontinousnegative}, we have
		\begin{align}
			\left[\left(1-\frac{C}{r-1}\right)\frac{t^2}{2r^2}+\frac{t}{2r}\right]\norm{A(X)}^2& \leqslant\E(t)\leqslant\E(0)=\frac{r-1}{2C}\norm{x_0-x^*}^2, \label{eq:divide1}\\
		\frac{t}{r}\inner{A(X),X-x^*}  &\leqslant\E(t)\leqslant\E(0)=\frac{r-1}{2C}\norm{x_0-x^*}^2.\label{eq:divide2}
		\end{align}
	By dividing $\left(1-\dfrac{C}{r-1}\right)\dfrac{t^2}{2r^2}+\dfrac{t}{2r}$ on both sides of \eqref{eq:divide1} and $\dfrac{t}{r}$ on both sides of \eqref{eq:divide2}, and then taking infimum respect to all $x^*\in\zero(A)$ on \eqref{eq:divide1}, we obtain \eqref{eq:t2bound} and \eqref{eq:t1bound}.
	
	Next, by the inequality \[
	\EE(t)\leqslant-\left(\frac{r-1-C}{r^2}t+1-\frac{1}{2r}\right)\norm{A(X)}^2-\frac{r-1}{r}\inner{A(X),X-x^*},
	\]
	we have
	\[\begin{gathered}
		\int_{0}^{\infty}\left(\frac{r-1-C}{r^2}t+1-\frac{1}{2r}\right)\norm{A(X(t))}^2\dd t \leqslant- \int_{0}^{\infty}\EE(t)\dd t\leqslant \E(0),\\
		\int_{0}^{\infty}\frac{r-1}{r}\inner{A(X(t)),X(t)-x^*}\dd t\leqslant-\int_{0}^{\infty}\EE(t)\dd t\leqslant\E(0).
	\end{gathered}
	\]
	The second inequality is \eqref{eq:integralfinite2}. Also, by taking infimum respect to all $x^*\in\zero(A)$ on first above inequality, we obtain \eqref{eq:integralfinite}.
	
\end{proof}
Next, we demonstrate that under the condition \( C < r - 1 \), the solution trajectory to \eqref{eq:acceleratedode} exhibits an \( o(1/t^2) \) convergence rate. Inspired by the analysis in \cite{yi23}, we first examine the properties of the functions \( X(t) - Z(t) \) and \( Z(t)-x^* \). Drawing on the proof technique from \cite{bot232} and the preceding arguments, we introduce the following pair of auxiliary functions to analyze \( X(t) - Z(t) \) and \( Z(t)-x^* \):
\begin{align}
	\G_1(t) &= \frac{1}{2}\norm{X(t) - Z(t)}^2 + \frac{(r - C)t + r^2}{2Ct}\norm{Z(t) - x^*}^2, \quad \forall x^* \in \zero(A), \label{eq:g1t}\\
	\G_2(t) &= \frac{1}{2}\norm{X(t) - Z(t)}^2 + \frac{t[(r - C)t + r^2]}{2r^2(1 + C)}\|\dot{X}(t)\|^2.\label{eq:g2t}
\end{align}
With these two functions, we can obtain the following result.

\begin{lemma}
	\label{lem:auxiliarycontinuous}
	Let \( A \) be a single-valued maximally monotone operator, and let \((X(t), Z(t))\) be the solution to \eqref{eq:acceleratedode}. If \( r > 1 \) and \( 0 < C < r - 1 \), then the following limits hold:
	\begin{align}
		&\lim_{t \to \infty} \norm{X(t) - Z(t)}^2 = 0, \label{eq:lemcontinuous1} \\
		&\lim_{t \to \infty} t^2 \|\dot{X}(t)\|^2 = 0, \label{eq:lemcontinuous2} \\
		&\lim_{t \to \infty} \norm{Z(t) - x^*}^2 \text{ exists for all } x^* \in \zero(A). \label{eq:lemcontinuous3}
	\end{align}
\end{lemma}
\begin{proof}
	First, we consider the derivative of $\G_2$. Directly calculating shows that:
	\begin{align*}
		\overset{\cdot}{\G}_2(t)= &\ \inner{\X-\Z,X-Z}+\frac{2(r-C)t+r^2}{2r^2(1+C)}\|\X\|^2+\frac{t[(r-C)t+r^2]}{r^2(1+C)}\inner{\X,\XX}\\
		=&\ \inner{\frac{r}{t}(Z-X)-\left(1+\frac{r}{t}-\frac{C}{r}\right)A(X),X-Z}+\frac{2(r-C)t+r^2}{2r^2(1+C)}\|\X\|^2\\
		&\ +\frac{(r-C)t+r^2}{r(1+C)}\inner{\X,-\frac{1+r}{r}\X-\left(1+\frac{t}{r}\right)\frac{\dd}{\dd t}A(X)-\frac{1+C}{r}A(X)}\\
		=&\ -\frac{r}{t}\norm{X-Z}^2-\left(1+\frac{r}{t}-\frac{C}{r}\right)\inner{A(X),X-Z}-\frac{2r(r-C)t+2r^3+r^2}{2r^2(1+C)}\|\X\|^2\\
		&\ -\frac{(r-C)t+r^2}{r^2}\inner{A(X),\X}-\frac{[(r-C)t+r^2](t+r)}{r^2(1+C)}\inner{\X,\frac{\dd}{\dd t}A(X)}\\
		=&\ -\frac{r}{t}\norm{X-Z}^2-\frac{2r(r-C)t+2r^3+r^2}{2r^2(1+C)}\|\X\|^2\\
		&\ -\frac{[(r-C)t+r^2](t+r)}{r^2(1+C)}\inner{\X,\frac{\dd}{\dd t}A(X)}+\frac{[(r-C)t+r^2](t+r)}{r^2t}\norm{A(X)}^2
	\end{align*}
	Let $\varepsilon>0$ be an arbitrary positive constant. By \eqref{eq:integralfinite}, we have $\int_{\varepsilon}^{\infty}\dfrac{[(r-C)u+r^2](u+r)}{r^2u}\norm{A(X(u))}^2\dd u<\infty$. Thus\[
	\frac{\dd}{\dd t}\left(\G_2(t)+\int_{t}^{\infty}\dfrac{[(r-C)u+r^2](u+r)}{r^2u}\norm{A(X(u))}^2\dd u\right)\leqslant 0,\quad\forall t>\varepsilon.
	\]
	In conclusion, $\lim\limits_{t\to\infty}\G_2(t)$ exists. Also, by integrating $\overset{\cdot}{\G}_2(t)$ from $\varepsilon$ to $\infty$, we have
	\begin{align*}
		&\ \int_{\varepsilon}^{\infty}\frac{r}{t}\norm{X(t)-Z(t)}^2+\frac{2(r-C)t+r^3+2r^2}{2r^3(1+C)}\|\X(t)\|^2\dd t\\
		\leqslant&\ \G_2(\varepsilon)+\int_{\varepsilon}^{\infty}\frac{[(r-C)t+r^2](t+r)}{r^2t}\norm{A(X(t))}^2\dd t<\infty.
	\end{align*}
	If $\lim\limits_{t\to\infty}\G_2(t)>0$, 
	\begin{align*}
		&\ \int_{\varepsilon}^{\infty}\frac{r}{t}\norm{X(t)-Z(t)}^2+\frac{2r(r-C)t+2r^3+r^2}{2r^2(1+C)}\|\X(t)\|^2\dd t\\
		\geqslant&\ \int_\varepsilon^\infty\frac{2r}{t}\G_2(t)\dd t=\infty,
	\end{align*}
	which leads to a contradiction. Therefore, $\lim\limits_{t\to\infty}\G_2(t)=0$, which implies\[
	\lim_{t\to\infty}t^2\|\X(t)\|^2 = 0,\quad\lim_{t\to\infty}\norm{X(t)-Z(t)}^2=0.
	\]
	Next, we show that $\lim\limits_{t\to\infty}\G_1(t)$ exists. The derivative of $\G_1(t)$ can be calculated as follows:\begin{align*}
		\frac{\dd}{\dd t}\G_1(t)=& \inner{\X-\Z, X-Z}+\frac{(r-C)t+r^2}{Ct}\inner{\Z, Z-x^*}-\frac{r^2}{2Ct^2}\norm{Z-x^*}^2  \\
		=& \inner{\frac{r}{t}(Z-X)-\left(1+\frac{r}{t}-\frac{C}{r}\right)A(X), X-Z}-\frac{r}{Ct^2}\norm{Z-x^*}^2\\
		&\ -\left(1+\frac{r}{t}-\frac{C}{r}\right)\inner{A(X), Z-x^*}\\
		=&\ -\frac{r}{t}\norm{X-Z}^2-\left(1+\frac{r}{t}-\frac{C}{r}\right)\inner{A(X),X-x^*}-\frac{r}{Ct^2}\norm{Z-x^*}^2\\
		\leqslant&\ 0.
	\end{align*}
	Due to non-negativity of $\G_1$,  $\lim\limits_{t\to\infty}\G_1(t)$ exists. Because of \eqref{eq:g1t}, the existence of $\lim\limits_{t\to\infty}\G_1(t)$ and \eqref{eq:lemcontinuous1}, $\lim\limits_{t\to\infty}\norm{Z(t)-x^*}^2$ exists for all $x^*\in\zero(A)$.
\end{proof}

By using Lemma \ref{lem:auxiliarycontinuous}, we can prove the $o(1/t^2)$ convergence rate of the solution trajectory to \eqref{eq:acceleratedode}.
\begin{theorem}[$o(1/t^2)$ convergence rates]
	Let $A$ be a single-valued maximally monotone operator, and let $(X(t),Z(t))$ be the solution to \eqref{eq:acceleratedode}. If $r>1$, $0<C< r-1$, then the following limits hold:\begin{align}
		\lim\limits_{t\to\infty}t^2\norm{A(X(t))}^2&=0,\label{eq:o1/t^2}\\
		\lim_{t\to\infty}t\inner{A(X(t)), X(t)-x^*}&=0,\quad\forall x^*\in\zero(A).
	\end{align}
\end{theorem}

\begin{proof}
	Theorem \ref{thm:convergenceratecontinuous} shows that $\E(t)$ is non-negative and non-increasing, which implies that \(\lim\limits_{t \to \infty} \E(t)\) exists. Next, by Theorem \ref{thm:convergenceratecontinuous} and the assumption \(0 < C < r - 1\), we have
	\begin{equation}
		\label{eq:t2infty}
		t^2 \norm{A(X)} < \infty, \quad \forall t > 0.
	\end{equation}
	Due to \eqref{eq:lemcontinuous1}, we have
	\begin{equation}
		\label{eq:limitcontinuous0}
		\lim_{t \to \infty} \frac{t}{r} \inner{A(X), X - Z} = 0.
	\end{equation}
	Based on \eqref{eq:lemcontinuous3}, \eqref{eq:t2infty}, and \eqref{eq:limitcontinuous0}, we conclude that
	\[
	\lim_{t \to \infty} \frac{t(t + r)}{2r^2} \norm{A(X)}^2 \text{ exists.}
	\]
	Based on this result and the inequality \eqref{eq:integralfinite}, we obtain \eqref{eq:o1/t^2}.
	
	Next, owing to \eqref{eq:lemcontinuous3} and \eqref{eq:o1/t^2}, we have
	\[
	\lim_{t \to \infty} \frac{1}{2} \norm{\frac{t}{r} \sqrt{\frac{C}{r-1}} A(X) + \sqrt{\frac{r-1}{C}} (x^* - Z)}^2 \text{ exists.}
	\]
	Thus, we have
	\[
	\lim_{t \to \infty} \frac{t}{r} \inner{A(X), X - x^*} \text{ exists.}
	\]
	Also, by using \eqref{eq:integralfinite2}, we have
	\[
	\lim_{t \to \infty} t \inner{A(X), X - x^*} = 0.
	\]
\end{proof}

So far, our focus has been on the convergence rates of the solution trajectory to \eqref{eq:acceleratedode}. Moving forward, we shift our attention to the convergence property of the solution trajectory to \eqref{eq:acceleratedode}.

\begin{theorem}[Weak Convergence Property]
	\label{thm:weakconvergence}
	Let \( A \) be a single-valued maximally monotone operator, and let \((X(t), Z(t))\) be the solution to \eqref{eq:acceleratedode}. If \( r > 1 \) and \( 0 < C < r - 1 \), then \( X(t) \) and \( Z(t) \) converge weakly to the same point in \(\zero(A)\).
\end{theorem}

\begin{proof}
	As a result of \eqref{eq:lemcontinuous1}, \eqref{eq:lemcontinuous3}, and Proposition \ref{prop:continuousopial}, it suffices to show that every weak cluster point of \( Z(t) \) belongs to \(\zero(A)\). Let \( z_\infty \) be any weak cluster point of \( Z(t) \), and let \(\{Z(t_n)\}\) be an arbitrary sequence such that \( Z(t_n) \rightharpoonup z_\infty \).
	
	Since \(\lim\limits_{t \to \infty} \norm{X(t) - Z(t)}^2 = 0\), it follows that \( X(t_n) \rightharpoonup z_\infty \). Additionally, by \eqref{eq:t2bound}, we have \( A(X(t_n)) \to 0 \). By Proposition \ref{prop:nonexpansiveweak}, we have \( 0 \in A(z_\infty) \), which means \( z_\infty \in \zero(A) \).
\end{proof}

\section{Symplectic Proximal Point Algorithm}
In this section, we will prove the convergence rate and weak convergence property of Algorithm \ref{al:sppa}. Similar to the argument in Section \ref{sec:2}, we first study the last-iterative convergence rate and ergodic convergence rate of Algorithm \ref{al:sppa} using the Lyapunov function technique. The Lyapunov function is given by:
\begin{equation}
	\E(k) = \frac{k(k + r)}{2r^2} \norm{\tilde{x}_k - x_k}^2 + \frac{k}{r} \inner{\tilde{x}_k - x_k, x_k - z_k} + \frac{r - 1}{2C} \norm{z_k - x^*}^2,
	\label{eq:lyapunovdiscrete}
\end{equation}
where $x^*\in\zero(A)$ and \(\tilde{x}_0 = x_0\). The function \(\E(k)\) can be viewed as the discrete counterpart of \eqref{eq:lyapunovcontinuous}. For simplicity, we define \(\overline{A}(x_k)\) as follows:
\[
\overline{A}(x_k) = \begin{cases}
	0, & \text{if } k = 0; \\
	\tilde{x}_k - x_k, & \text{if } k \geq 1.
\end{cases}
\]
By the definition of the proximal point operator, we have \(\overline{A}(x_{k+1}) \in A(x_{k+1})\), \(\forall k \geqslant 0\).

First, we study the non-negativity and non-increasing property of \(\E(k)\).

\begin{lemma}
	\label{lem:lyapunovdiscretenonnegative}
	Let $A$ be a maximally monotone operator, and let $\{x_k\}$ and $\{z_k\}$ be the sequences generated by Algorithm \ref{al:sppa}. The Lyapunov function $\E(k)$ given in \eqref{eq:lyapunovdiscrete} can be represented as\begin{align*}
	\E(k)=&\ \frac{(r-1-C)k^2+kr(r-1)}{2r^2(r-1)}\norm{\overline{A}(x_k)}^2+\frac{k}{r}\inner{\overline{A}(x_k),x_k-x^*}\\
	&\ +\frac{1}{2}\norm{\frac{k}{r}\sqrt{\frac{C}{r-1}}\overline{A}(x_k)+\sqrt{\frac{r-1}{C}}(x^*-z_k)}^2.
	\end{align*}
	If $0<C\leqslant r-1$, $\E(k)$ can be represented by a sum of three non-negative terms.
\end{lemma}
\begin{proof}
	Similar to the proof for Lemma \ref{lem:lyapunovcontinousnegative}, we have that\begin{align*}
		\E(k)= &\ \frac{k(k+r)}{2r^2}\norm{\overline{A}(x_k)}^2+\frac{k}{r}\inner{\overline{A}(x_k),x_k-x^*} \\
		 &\ +\frac{k}{r}\inner{\overline{A}(x_k),x^*-z_k}+\frac{r-1}{2C}\norm{z_k-x^*}^2\\
		 =&\ \frac{(r-1-C)k^2+kr(r-1)}{2r^2(r-1)}\norm{\overline{A}(x_k)}^2+\frac{k}{r}\inner{\overline{A}(x_k),x_k-x^*}\\
		 &\ +\frac{1}{2}\norm{\sqrt{\frac{Ck^2}{r^2(r-1)}}\overline{A}(x_k)+\sqrt{\frac{r-1}{C}}(x^*-z_k)}^2.
	\end{align*}
	Since $0<C\leqslant r-1$, we have\[
	\frac{(r-1-C)k^2+kr(r-1)}{2r^2(r-1)}\geqslant 0,\quad\forall k\geqslant 0.
	\]
	Since $A$ is monotone, $\E(k)$ can be represented by a sum of three non-negative terms.
\end{proof}

\begin{lemma}
	\label{lem:lyapunovdiscrete}
	Let $A$ be a maximally monotone operator, and let $\{x_k\}$ and $\{z_k\}$ be the sequences generated by Algorithm \ref{al:sppa}. If $r>1$ and $C\leqslant r-1$, then we have\begin{equation}
		\label{eq:differenceestimation}
		\begin{split}
			&\ \E(k+1)-\E(k)\\
			\leqslant&\ -\frac{r-1}{r}\inner{\overline{A}(x_{k+1}),x_{k+1}-x^*}-\frac{(2r-2-2C)k+r^2+r-2}{2r^2}\norm{\overline{A}(x_{k+1})}^2\\
			\leqslant&\ 0.
		\end{split}
	\end{equation}
\end{lemma}
\begin{proof}
	Step 1: Dividing the difference $\E(k+1)-\E(k)$ into three parts.\begin{align*}
		\E(k+1)-\E(k)= &\  \frac{(k+1)(k+1+r)}{2r^2}\norm{\overline{A}(x_{k+1})}^2-\frac{k(k+r)}{2r^2}\norm{\overline{A}(x_k)}^2\\
		&\ +\underbrace{\frac{k+1}{r}\inner{\overline{A}(x_{k+1}),x_{k+1}-z_{k+1}}-\frac{k}{r}\inner{\overline{A}(x_k),x_k-z_k}}_{\text{I}}\\
		&\ +\underbrace{\frac{r-1}{2C}\left(\norm{z_{k+1}-x^*}^2-\norm{z_k-x^*}^2\right)}_{\text{II}}.
	\end{align*}
	
	Step 2: Transforming I and II. First, we focus on I. \begin{align*}
		\text{I}= &\ \frac{1}{r}\inner{\overline{A}(x_{k+1}),x_{k+1}-z_{k+1}}+\frac{k}{r}\inner{\overline{A}(x_{k+1}), x_{k+1}-x_k-z_{k+1}+z_k} \\
		&\ +\frac{k}{r}\inner{\overline{A}(x_{k+1})-\overline{A}(x_k),x_k-z_k}.
	\end{align*}
	First of all, since $z_{k+1}=z_k-\dfrac{C}{r}\overline{A}(x_{k+1})$, the first term of I can be replaced by \[
	\frac{1}{r}\inner{\overline{A}(x_{k+1}),x_{k+1}-z_{k+1}}=\frac{1}{r}\inner{\overline{A}(x_{k+1}),x_{k+1}-z_k}+\frac{C}{r^2}\norm{\overline{A}(x_{k+1})}^2.
	\] 
	By using the equality \[
	k(x_{k+1}-x_k)=r(z_k-x_{k+1})-(k+r)\overline{A}(x_{k+1}),
	\]we have\begin{align*}
		&\ \frac{k}{r}\inner{\overline{A}(x_{k+1}), x_{k+1}-x_k-z_{k+1}+z_k}\\
		=&\ \frac{1}{r}\inner{\overline{A}(x_{k+1}),-\left(k+r-\frac{kC}{r}\right)\overline{A}(x_{k+1})+r(z_k-x_{k+1})}.	 
	\end{align*}
	Also, by the equality \[
	z_k=r^{-1}\left[(k+r)\overline{A}(x_{k+1})+(k+r)x_{k+1}-kx_k\right],
	\]we have
	\begin{align*}
		&\ \frac{k}{r}\inner{\overline{A}(x_{k+1})-\overline{A}(x_k),x_k-z_k}\\
		=&\ \frac{k}{r}\inner{\overline{A}(x_{k+1})-\overline{A}(x_k), -\left(1+\frac{k}{r}\right)(x_{k+1}-x_k)-\left(1+\frac{k}{r}\right)\overline{A}(x_{k+1})}. 
	\end{align*}
	Next, by using the equality\[
	\inner{a,b}=\frac{1}{2}\norm{a}^2+\frac{1}{2}\norm{b}^2-\frac{1}{2}\norm{a-b}^2,\quad\forall a, b\in\hh,
	\]
	we can transform I into the following form:\begin{align*}
		\text{I}= &\  -\frac{r-1}{r}\inner{\overline{A}(x_{k+1}), x_{k+1}-z_k}-\frac{k(k+r)}{r^2}\inner{\overline{A}(x_{k+1})-\overline{A}(x_k),x_{k+1}-x_k}\\
		&\ -\left[\frac{k}{r}+1+\frac{k(k+r)}{2r^2}-\frac{kC}{r^2}\right]\norm{\overline{A}(x_{k+1})}^2\\
		&\ +\frac{k(k+r)}{2r^2}\left(\norm{\overline{A}(x_k)}^2-\norm{\overline{A}(x_{k+1})-\overline{A}(x_k)}^2\right).
	\end{align*}
	
	Now turning to consider II. \begin{align*}
		\text{II} & =\frac{r-1}{C}\inner{z_{k+1}-z_k,z_k-x^*}+\frac{r-1}{2C}\norm{z_{k+1}-z_k}^2\\
		&=-\frac{r-1}{r}\inner{\overline{A}(x_{k+1}),z_k-x^*}+\frac{(r-1)C}{2r^2}\norm{\overline{A}(x_{k+1})}^2.
	\end{align*}
	
	Step 3: Estimating the upper bound of $\E(k+1)-\E(k)$. As a consequence of previous calculation, we have\begin{align*}
		&\ \E(k+1)-\E(k)\\
		= &\ -\frac{r-1}{r}\inner{\overline{A}(x_{k+1}),x_{k+1}-x^*}-\frac{k(k+r)}{r^2}\inner{\overline{A}(x_{k+1})-\overline{A}(x_k),x_{k+1}-x_k}\\
		&\ -\frac{(2r-2-2C)k+2r^2-r-1-C(r-1)}{2r^2}\norm{\overline{A}(x_{k+1})}^2\\
		&\ -\frac{k(k+r)}{2r^2}\norm{\overline{A}(x_{k+1})-\overline{A}(x_k)}^2.
	\end{align*}
	By the assumptions that $r>1$ and $C\leqslant r-1$, we have $2r-2-2C\geqslant 0$ and $2r^2-r-1-C(r-1)\geqslant 2r^2-r-1-(r-1)^2=r^2+r-2>0$. Since $A$ is monotone and $x^*\in\zero(A)$, we have\[
	\inner{\overline{A}(x_{k+1}),x_{k+1}-x^*}\geqslant 0,\quad\inner{\overline{A}(x_{k+1})-\overline{A}(x_k),x_{k+1}-x_k}\geqslant 0.
	\]
	Thus we obtain \eqref{eq:differenceestimation}.
\end{proof}

With previous analysis of $\E(k)$, we can obtain the following theorem.
\begin{theorem}
	\label{thm:symplecticpparate}
	Let $A$ be a maximally monotone operator, and let $\{x_k\}$ and $\{z_k\}$ be the sequences generated by Algorithm \ref{al:sppa}. If $r>1$ and $0<C\leqslant r-1$, then $\{\E(k)\}$ is non-increasing and non-negative. The last-iterate convergence rates of Algorithm \ref{al:sppa} are given as follows:
	\begin{align}
		\norm{\overline{A}(x_k)}^2&\leqslant\frac{r^2(r-1)^2}{[C(r-1)-C^2]k^2+Cr(r-1)k}\cdot\dist(x_0,\zero(A))^2,\label{eq:k2bound}\\
		\inner{\overline{A}(x_k),x_k-x^*}&\leqslant\frac{r^2-r}{2Ck}\cdot\norm{x_0-x^*}^2,\quad\forall x^*\in\zero(A).\label{eq:k1bound}
	\end{align}
	In addition, we have
		\begin{align}
			\sum_{k=0}^\infty\frac{(2r-2-2C)k+r^2-1}{2r^2}\norm{\overline{A}(x_{k+1})}^2\leqslant \frac{r-1}{2C}\cdot\dist(x_0,\zero(A))^2.\label{eq:ksumbound}\\
		 \sum_{k=0}^{\infty}\frac{r-1}{r}\inner{\overline{A}(x_k),x_{k+1}-x^*}\leqslant\frac{r-1}{2C}\cdot\norm{x_0-x^*}^2,\quad\forall x^*\in\zero(A),\label{eq:sumbound}
		\end{align}
\end{theorem}
\begin{proof}
	Due to Lemma \ref{lem:lyapunovdiscrete}, we have $\{\E(k)\}$ is non-increasing. Then by Lemma \ref{lem:lyapunovdiscretenonnegative}, we have\begin{align}
		\frac{(r-1-C)k^2+kr(r-1)}{2r^2(r-1)}\norm{\overline{A}(x_k)}^2 & \leqslant\E(k)\leqslant\E(0)=\frac{r-1}{2C}\norm{x_0-x^*}^2,\label{eq:divide3} \\
		\frac{k}{r}\inner{\overline{A}(x_k),x_k-x^*} & \leqslant\E(k)\leqslant\E(0)=\frac{r-1}{2C}\norm{x_0-x^*}^2.\label{eq:divide4}
	\end{align}
	By dividing $\dfrac{(r-1-C)k^2+kr(r-1)}{2r^2(r-1)}$ on both sides of \eqref{eq:divide3}, $\dfrac{k}{r}$ on both sides of  \eqref{eq:divide4} and taking infimum respect to all $x^*\in\zero(A)$ on \eqref{eq:divide3}, we obtain \eqref{eq:k2bound} and \eqref{eq:k1bound}. Also, due to \eqref{eq:differenceestimation}, we have\begin{align*}
		\sum_{k=0}^{\infty}\frac{(2r-2-2C)k+r^2-1}{2r^2}\norm{\overline{A}(x_{k+1})}^2&\leqslant\sum_{k=0}^{\infty}\E(k)-\E(k+1)\\
		&\leqslant\E(0)=\frac{r-1}{2C}\norm{x_0-x^*}^2,\\
		\sum_{k=0}^{\infty}\frac{r-1}{r}\inner{\overline{A}(x_{k+1}),x_{k+1}-x^*} & \leqslant\sum_{k=0}^{\infty}\E(k)-\E(k+1)\\
		&\leqslant\E(0)=\frac{r-1}{2C}\norm{x_0-x^*}^2. \\
	\end{align*}
	Taking infimum respect to all $x^*\in\zero(A)$ on first above inequality, we obtain \eqref{eq:ksumbound} and \eqref{eq:sumbound}.
\end{proof}

Next, we show that the \( o(1/k^2) \) convergence rates and weak convergence property of Algorithm \ref{al:sppa} hold when \( C < r - 1 \). Similar to the auxiliary functions introduced in Section \ref{sec:2}, we consider the following two auxiliary sequences:
\begin{align}
	\G_1(k) & = \frac{1}{2}\norm{x_k-z_k}^2+\frac{(r-C)(k+r)^2}{k^2C}\norm{z_k-z^*}^2,\label{eq:g1k} \\
	\G_2(k) & = \frac{1}{2}\norm{x_k-z_k}^2+\frac{(r-C)(k-1)^2(k+r)^2}{2(1+C)k^2r^2}\norm{x_k-x_{k-1}}^2,\quad k\geqslant 1.\label{eq:g2k}
\end{align}
\begin{lemma}
	\label{lem:auxiliarydiscrete}
	Let $\{x_k\}$ and $\{z_k\}$ be the sequences generated by Algorithm \ref{al:sppa}. If $r>1$, $0<C<r-1$, we have\begin{align}
		&\lim_{k\to\infty}\norm{x_k-z_k}^2=0,\label{eq:lemdiscrete1}\\
		&\lim_{k\to\infty}k^2\norm{x_{k+1}-x_k}^2=0,\label{eq:lemdiscrete2}\\
		&\lim_{k\to\infty}\norm{z_k-x^*}^2\text{ exists for all }x^*\in\zero(A).\label{eq:lemdiscrete3}
	\end{align}
\end{lemma}
\begin{proof}
	First, we consider the differences $\dfrac{1}{2}\norm{x_{k+1}-z_{k+1}}^2-\dfrac{1}{2}\norm{x_k-z_k}^2$ and  $k^2\norm{x_{k+1}-x_k}^2-(k-1)^2\norm{x_k-x_{k-1}}^2$. Direct calculation shows that:
	\begin{align*}
		&\ \dfrac{1}{2}\norm{x_{k+1}-z_{k+1}}^2-\dfrac{1}{2}\norm{x_k-z_k}^2\\
		=&\ \frac{1}{2}\inner{x_{k+1}+x_k-z_{k+1}-z_k, x_{k+1}-x_k-z_{k+1}+z_k}\\
		=&\ -\frac{(2k+r)r}{2k^2}\norm{x_{k+1}-z_{k+1}}^2-\frac{(r^2-C^2)(k+r)^2}{2k^2r^2}\norm{\overline{A}(x_{k+1})}^2\\
		&\ -\frac{(r-C)(k+r)^2}{k^2r}\inner{\overline{A}(x_{k+1}),x_{k+1}-z_k},\\
		&\frac{k^2}{2}\norm{x_{k+1}-x_k}^2-\frac{(k-1)^2}{2}\norm{x_k-x_{k-1}}^2 \\
		\leqslant&\ k\inner{x_{k+1}-x_k, k(x_{k+1}-x_k)-(k-1)(x_k-x_{k-1})}\\
		=&\ k\inner{x_{k+1}-x_k, -(k+r)\overline{A}(x_{k+1})-r(x_{k+1}-x_k)+(k+r-1-C)\overline{A}(x_k)}\\
		=&\ -(k+r-1-C)k\inner{x_{k+1}-x_k, \overline{A}(x_{k+1})-\overline{A}(x_k)}-kr\norm{x_{k+1}-x_k}^2\\
		&\ -(1+C)k\inner{\overline{A}(x_{k+1}),x_{k+1}-x_k}\\
		=&\ -(k+r-1-C)k\inner{x_{k+1}-x_k, \overline{A}(x_{k+1})-\overline{A}(x_k)}-kr\norm{x_{k+1}-x_k}^2\\
		&\ +(1+C)r\inner{\overline{A}(x_{k+1}),x_{k+1}-z_k}+(1+C)(k+r)\norm{\overline{A}(x_{k+1})}^2\\
		\leqslant&\ -kr\norm{x_{k+1}-x_k}^2+(1+C)r\inner{\overline{A}(x_{k+1}),x_{k+1}-z_k}\\
		&\ +(1+C)(k+r)\norm{\overline{A}(x_{k+1})}^2.
	\end{align*}
	Based on previous calculation and the inequality $\dfrac{(k+1+r)^2}{(k+1)^2}<\dfrac{(k+r)^2}{k^2}$, we have\begin{align*}
		&\ \G_2(k+1)-\G_2(k)\\
		\leqslant&\ \frac{1}{2}\norm{x_{k+1}-z_{k+1}}^2-\frac{1}{2}\norm{x_k-z_k}^2\\
		&\ +\frac{(r-C)(k+r)^2}{(1+C)k^2r^2}\left(\frac{k^2}{2}\norm{x_{k+1}-x_k}^2-\frac{(k-1)^2}{2}\norm{x_k-x_{k-1}}^2\right)\\
		\leqslant&\ -\frac{(2k+r)r}{2k^2}\norm{x_{k+1}-z_{k+1}}^2-\frac{(r^2-C^2)(k+r)^2}{2k^2r^2}\norm{\overline{A}(x_{k+1})}^2\\
		&\ -\frac{(r-C)(k+r)^2}{(1+C)kr}\norm{x_{k+1}-x_k}^2+\frac{(r-C)(k+r)^3}{k^2r^2}\norm{\overline{A}(x_{k+1})}^2.
	\end{align*}
	By \eqref{eq:ksumbound}, we have\[
	\sum_{k=1}^{\infty}\frac{(r-C)(k+r)^3}{k^2r^2}\norm{\overline{A}(x_{k+1})}^2<\infty.
	\]
	Consequently, we have
	\begin{align*}
		&\ \G_2(k+1)+\sum_{i=k+1}^\infty\frac{(r-C)(i+r)^3}{i^2r^2}\norm{\overline{A}(x_{i+1})}^2\\
		\leqslant&\  \G_2(k)+\sum_{i=k}^\infty\frac{(r-C)(i+r)^3}{i^2r^2}\norm{\overline{A}(x_{i+1})}^2.
	\end{align*}
	Thus $\lim\limits_{k\to\infty}\mathcal{G}_2(k)$ exist. Also, we have
	\begin{align*}
		&\ \frac{(2k+r)r}{2k^2}\norm{x_{k+1}-z_{k+1}}^2+\frac{(r-C)(k+r)^2}{(1+C)kr}\norm{x_{k+1}-x_k}^2\\
		\leqslant&\ \G_2(k)-\G_2(k+1)+\frac{(r-C)(k+r)^3}{k^2r^2}\norm{\overline{A}(x_{k+1})}^2.
	\end{align*}
	Then we have
	\begin{align*}
		&\ \sum_{k=1}^\infty\frac{(2k+r)r}{2k^2}\norm{x_{k+1}-z_{k+1}}^2+\frac{(r-C)(k+r)^2}{(1+C)kr}\norm{x_{k+1}-x_k}^2\\
		\leqslant&\ \G_2(1)+\sum_{k=1}^\infty\frac{(r-C)(k+r)^3}{k^2r^2}\norm{\overline{A}(x_{k+1})}^2<\infty.
	\end{align*}
	Moreover, based on above estimation, the inequality \[
	\frac{(2k+r)r}{2k^2}\norm{x_{k+1}-z_{k+1}}^2+\frac{(r-C)(k+r)^2}{(1+C)kr}\norm{x_{k+1}-x_k}^2\geqslant 2r\G_2(k+1),
	\]
	and the existence of $\lim\limits_{k\to\infty}\G_2(k)$, we have
	$\lim\limits_{k\to\infty}\G_2(k)=0$, 
	which implies that\[
	\lim_{k\to\infty}\norm{x_k-z_k}^2=0,\quad\lim_{k\to\infty}k^2\norm{x_{k+1}-x_k}^2=0.
	\]
	
	Next, we show that $\lim\limits_{k\to\infty}\norm{z_k-x^*}^2$ exists for all $x^*\in\zero(A)$. First, we show that $\{\G_1(k)\}$ is non-increasing. \begin{align*}
		&\ \G_1(k+1)-\G_1(k)\\
		\leqslant &\  -\frac{(2k+r)r}{2k^2}\norm{x_{k+1}-z_{k+1}}^2-\frac{(r^2-C^2)(k+r)^2}{2k^2r^2}\norm{\overline{A}(x_{k+1})}^2\\
		&\ -\frac{(r-C)(k+r)^2}{k^2r}\inner{\overline{A}(x_{k+1}),x_{k+1}-z_k}\\
		&\ +\frac{(r-C)(k+r)^2}{2k^2C}\left(\norm{z_{k+1}-x^*}^2-\norm{z_k-z^*}^2\right)\\
		\leqslant&\ -\frac{(2k+r)r}{2k^2}\norm{x_{k+1}-z_{k+1}}^2-\frac{(r-C)(k+r)^2}{2k^2r}\norm{\overline{A}(x_{k+1})}^2\\
		&\ -\frac{(r-C)(k+r)^2}{k^2r}\inner{\overline{A}(x_{k+1}),x_{k+1}-x^*}\\
		\leqslant&\ 0.
	\end{align*}
	Owing to non-negativity of $\G_1(k)$, $\lim\limits_{k\to\infty}\mathcal{G}_1(k)$ exists. Since $\lim\limits_{k\to\infty}\norm{x_k-z_k}^2=0$,  $\lim\limits_{k\to\infty}\norm{z_k-x^*}^2$ exists for all $x^*\in\zero(A)$.
\end{proof}
\begin{theorem}[$o(1/k^2)$ convergence rate] 
	Let $A: \hh\to 2^{\hh}$ be a maximally monotone operator, and let $\{x_k\}$, $\{z_k\}$, $\{\tilde{x}_{k+1}\}$ be the sequences generated by Algorithm \ref{al:sppa}. If $r>1$ and $0<C<r-1$, we have\begin{align}
		\lim_{k\to\infty}k^2\norm{\overline{A}(x_k)}^2 &=0,\\
		\lim_{k\to\infty}k\inner{\overline{A}(x_k), x_k-x^*} &=0,\quad\forall x^*\in\zero(A).
	\end{align}
\end{theorem}
\begin{proof}
	Because of Theorem \ref{thm:symplecticpparate},  $\lim\limits_{k\to\infty}\E(k)$ exists. Next, by \eqref{eq:k2bound} and \eqref{eq:lemdiscrete1}, we have\[
	\lim_{k\to\infty}\frac{k}{r}\inner{\overline{A}(x_k),x_k-z_k}=0.
	\]
	Also, by using the existence of $\lim\limits_{k\to\infty}\E(k)$ and \eqref{eq:lemdiscrete3}, we have\[
	\lim_{k\to\infty}k^2\norm{\overline{A}(x_k)}^2\text{ exists}.
	\]
	Additionally, owing to \eqref{eq:ksumbound}, we have\[
	\lim\limits_{k\to\infty}k^2\norm{\overline{A}(x_k)}^2=0.
	\]
	Since $\lim\limits_{k\to\infty}\norm{z_k-x^*}^2$ exists for all $x^*\in\zero(A)$, $\{z_k\}$ is bounded. By \eqref{eq:lemdiscrete1}, $\{x_k\}$ is bounded, too. Thus \[
	\lim_{k\to\infty}\frac{k}{r}\inner{\overline{A}(x_k),x_k-x^*}=0.
	\]
\end{proof}

\begin{theorem}[Weak Convergence Property]
	\label{thm:sppaweakconvergence}
	Let \( A \) be a maximally monotone operator, and let \(\{x_k\}\), \(\{z_k\}\),  \(\{\tilde{x}_{k+1}\}\) be the sequences generated by Algorithm \ref{al:sppa}. If \( r > 1 \) and \( 0 < C < r - 1 \), then both sequences \(\{x_k\}\) and \(\{z_k\}\) converge weakly to the same point in \(\zero(A)\).
\end{theorem}

\begin{proof}
	Owing to Lemma \ref{lem:auxiliarydiscrete} and Proposition \ref{prop:discreteopial}, it suffices to show that every weak cluster point of \(\{z_k\}\) belongs to \(\zero(A)\). Let \( z_\infty \) be any weak cluster point of \(\{z_k\}\), and let \(\{z_{k_j}\}\) be an arbitrary sequence of \(\{z_k\}\) such that \( z_{k_j} \rightharpoonup z_\infty \).
	
	Since \(\lim\limits_{k \to \infty} \norm{x_k - z_k}^2 = 0\), it follows that \( x_{k_j} \rightharpoonup z_\infty \). Given that \(\overline{A}(x_{k_j}) \in A(x_{k_j})\) and \(\overline{A}(x_{k_j}) \to 0\), we have \( 0 \in A(z_\infty) \) by Proposition \ref{prop:nonexpansiveweak}, which means \( z_\infty \in \zero(A) \).
\end{proof}

\section{Numerical Experiments}
In this section, we present several numerical experiments to illustrate the performance of our SPPA. Additionally, we demonstrate how the parameters $r$ and $C$ in SPPA influence the numerical performance.
\subsection{Matrix Example}
\label{sec:num1}

First, we consider the maximally monotone operator defined as
\begin{equation}
	A = \begin{pmatrix}
		\mathbb{O} & \mathbb{I}_d \\
		-\mathbb{I}_d & \mathbb{O}
	\end{pmatrix},
	\label{eq:matrixexample}
\end{equation}
where \(\mathbb{I}_d\) is the identity matrix and \(\mathbb{O}\) is the zero matrix in \(\mathbb{R}^{d \times d}\). This operator has been used in the literature to illustrate the worst-case performance of the proximal point method (see \cite{gu20,park22}). It is easy to verify that \(\zero(A) = \{0\}\).

In our experiment, we test the numerical performance of SPPA with \( r = 2 \) and \( C = 1 \). The accelerated PPA \eqref{eq:acceleratedproximal} and the Fast K-M algorithm with \( s = 2 \) and \(\alpha = 3 \) \eqref{eq:fastkm} are comparable. Here, \( d \) is set to 1000, and the initial point for all these methods is \(\begin{pmatrix}
	\vec{1}_d \\
	\vec{0}_d
\end{pmatrix}\), where $\vec{1}_d$ and $\vec{0}_d$ denote all ones vector and all zeros vector in $\mathbb{R}^d$, respectively. The PPA is not included in our experiment because its convergence rate in this case is \( O(\sqrt{2}^{-k}) \).

\begin{figure}[ht]
	\centering
	\subfigure{\includegraphics[width=.45\textwidth]{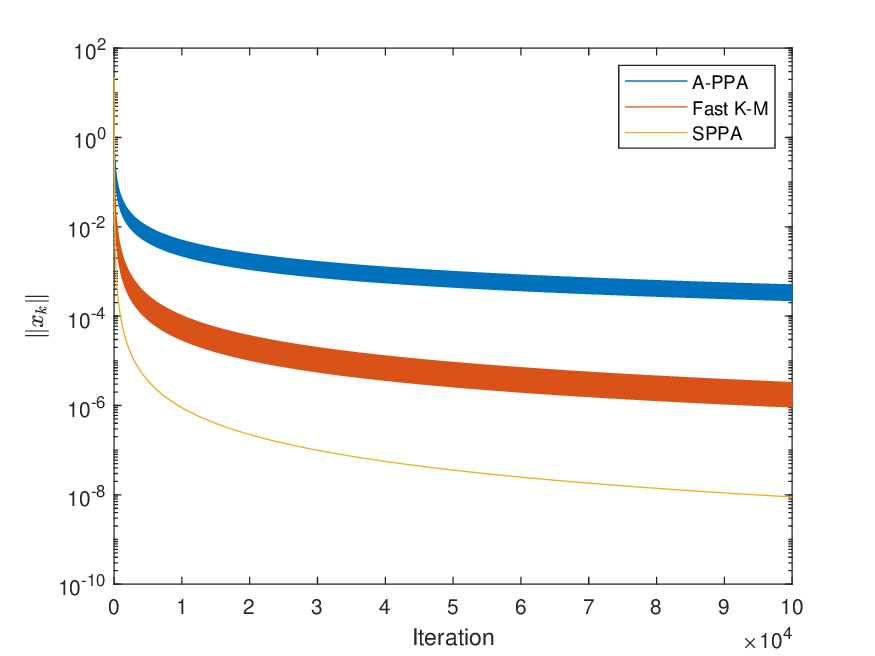}}
	\subfigure{\includegraphics[width=.45\textwidth]{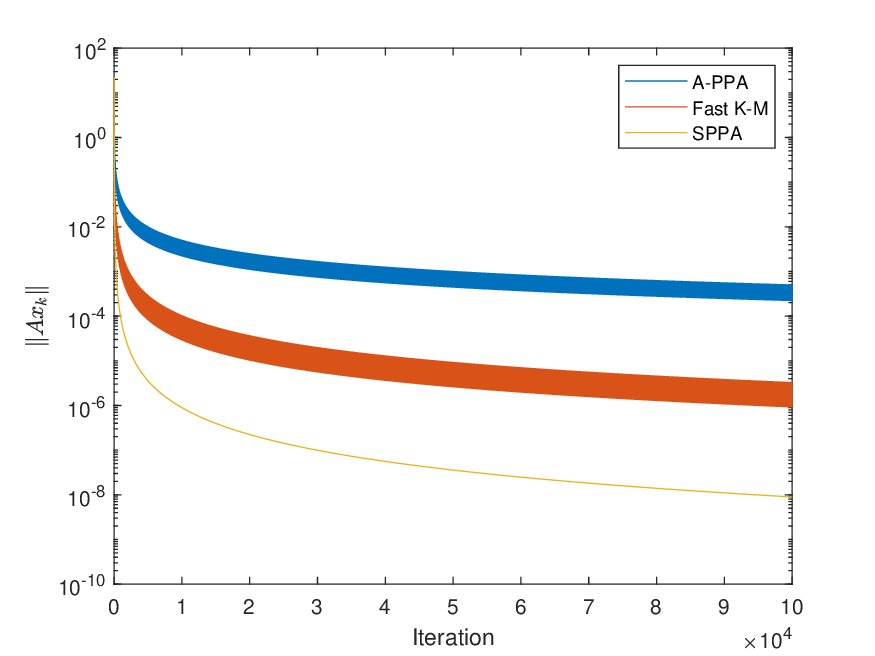}}
	\caption{Comparison between several accelerated PPAs.}
	\label{fig:num1}
\end{figure}

As we can see from Figure \ref{fig:num1}, our SPPA converges faster than accelerated PPA and the Fast K-M algorithm. In addition, the oscillatory behavior of accelerated PPA and the Fast K-M algorithm can be observed. By contrast, the SPPA does not suffer from oscillation, which may explain why it converges faster than the other two algorithms.

\begin{figure}[ht]
	\centering
	\subfigure{\includegraphics[width=.45\textwidth]{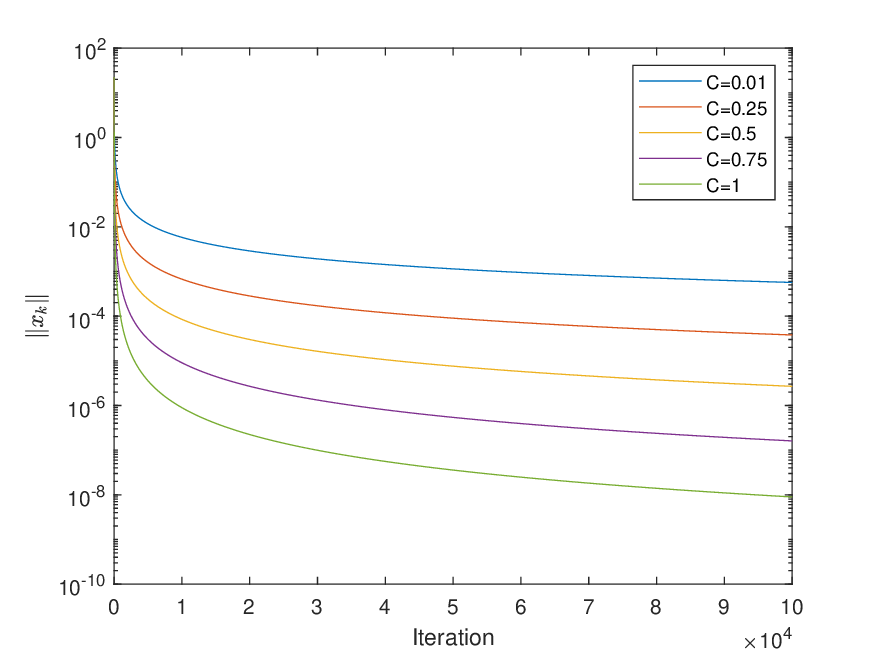}}
	\subfigure{\includegraphics[width=.45\textwidth]{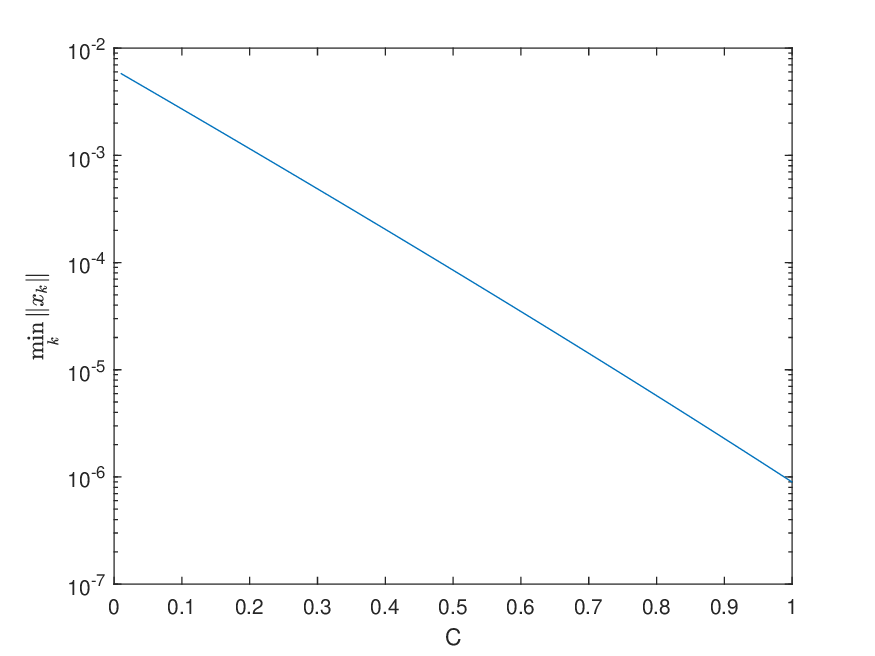}}
	\caption{Comparison between SPPAs with different $C$.}
	\label{fig:num1D}
\end{figure}

Next, we conduct a new experiment to show how the parameter \( C \) influences the numerical performance of SPPA. In this experiment, the maximally monotone operator is still given by \eqref{eq:matrixexample}, and the initial point is set to \(\begin{pmatrix}
	\vec{1}_d \\
	\vec{0}_d
\end{pmatrix}\). The parameter \( r \) is set to  2. We run SPPAs with \( C = 0.01 \), \( C = 0.25 \), \( C = 0.5 \), \( C = 0.75 \) and \( C = 1 \) for $10^6$ iterations, and we plot \(\norm{x_k}\) against the iteration number as shown in the left-hand side of Figure \ref{fig:num1D}. Additionally, we run SPPAs with different values of \( C \) for \( 10^4 \) iterations and plot \(\min_{1 \leqslant k \leqslant 10^4} \norm{x_k}\) against \( C \) as shown in the right-hand side of Figure \ref{fig:num1D}.

As we can see from Figure \ref{fig:num1D}, the numerical performance of SPPA improves as \( C \) increases. Fortunately, the result in Theorem \ref{thm:symplecticpparate} can help explain this phenomenon. Let us revisit the inequality \eqref{eq:k1bound}. The right-hand side of \eqref{eq:k1bound} decreases as \( C \) increases, suggesting that a larger \( C \) leads to faster convergence of SPPA. However, \eqref{eq:k2bound} indicates that \( C = 0.5(r-1) \) is optimal. How do we reconcile this apparent conflict?

One possible explanation is that the term \(\inner{\overline{A}(x_k), x_k - x^*}\), which appears on the left-hand side of \eqref{eq:k1bound}, characterizes the distance between \( x_k \) and \( x^* \). For instance, if \( A \) is the derivative of a convex function \(\varphi\), then \(\inner{\overline{A}(x_k), x_k - x^*}\) provides an upper bound on the Bergman divergence (see \cite{beck03, aharon2001} for an introduction to Bergman divergence) between \( x_k \) and \( x^* \). In practice, the distance between \( x_k \) and \( x^* \) is often more important.

\begin{figure}[ht]
	\centering
	\subfigure{\includegraphics[width=.45\textwidth]{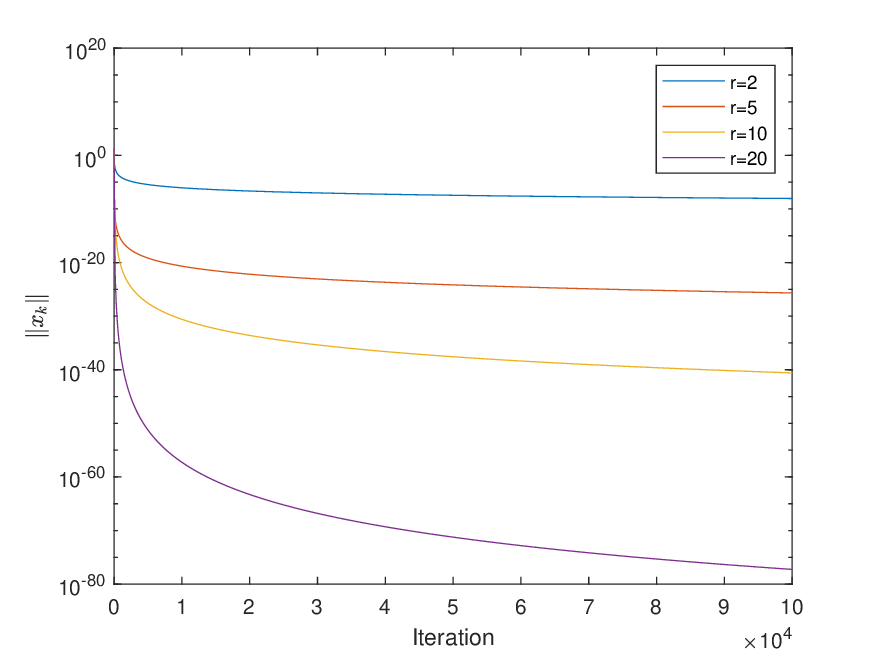}}
	\subfigure{\includegraphics[width=.45\textwidth]{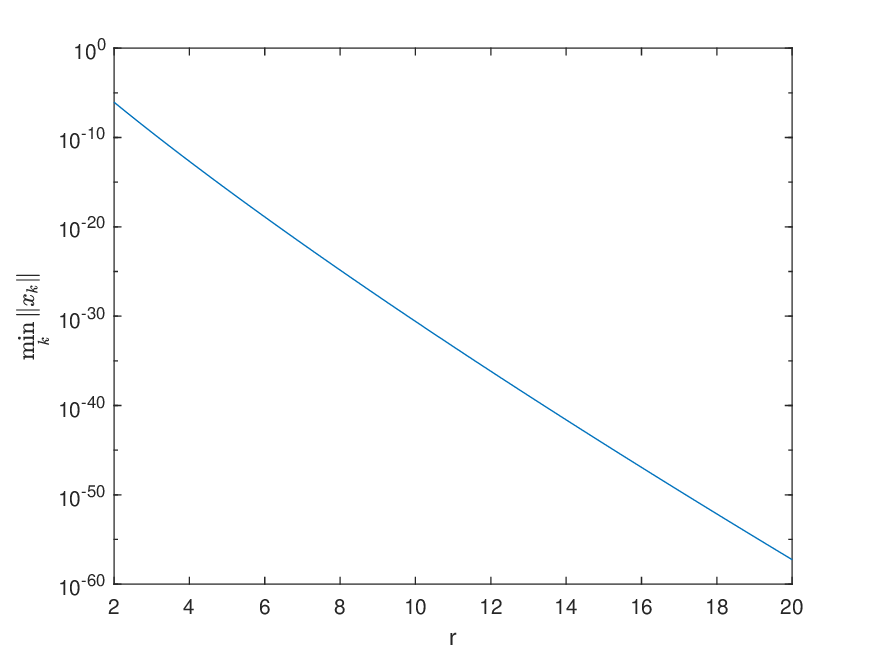}}
	\caption{Comparison between SPPAs with different $r$.}
	\label{fig:num1R}
\end{figure}

Subsequently, we present how \( r \) influences the numerical performance of SPPA. Similar to the previous experiment, we run SPPAs with different values of \( r \) and set \( C = r - 1 \), resulting in Figure \ref{fig:num1R}.

As we can see from Figure \ref{fig:num1R}, a larger \( r \) leads to faster convergence of SPPA. To explain this phenomenon, let us revisit the inequality \eqref{eq:sumbound}. Since \( C = r - 1 \), the right-hand side of \eqref{eq:sumbound} remains constant. Meanwhile, the coefficient of \(\inner{\overline{A}(x_k), x_k - x^*}\) on the left-hand side of \eqref{eq:sumbound} increases as \( r \) becomes larger. Therefore, a larger \( r \) leads to a faster ergodic convergence rate of \(\inner{\overline{A}(x_k), x_k - x^*}\), which ultimately results in faster convergence of SPPA.

\subsection{Intersection of Two Closed Convex Sets Example}
In this section, we consider the following problem:
\begin{equation}
	x \in C_1 \cap C_2,
	\label{eq:intersection}
\end{equation}
where \( C_1 \) and \( C_2 \) are two closed convex subsets. The parallel projection algorithm (Example 5.22 in \cite{bauschke17}) is an effective method to solve \eqref{eq:intersection}. The recursive rule of the parallel projection algorithm is given by:
\begin{equation}
	x_{k+1} = x_k + \lambda_k \left( \omega_1 P_{C_1}(x_k) + (1 - \omega_1) P_{C_2}(x_k) - x_k \right),
	\label{eq:parallelprojection}
\end{equation}
where \( P_{C_i} \) denotes the projection onto \( C_i \) for \( i = 1, 2 \), \( \omega_1 \in (0, 1) \), and \( \{\lambda_k\} \) is a sequence in \([0, 2]\) such that \( \sum_{k=0}^{\infty} \lambda_k (2 - \lambda_k) = \infty \). The operator \( T = \omega_1 P_{C_1} + (1 - \omega_1) P_{C_2} \) has been proven to be a proximal point operator for a specific maximally monotone operator. Thus, the parallel projection algorithm can be accelerated from the perspective of the accelerated PPA.

\begin{figure}[ht]
	\centering
	\subfigure{\includegraphics[width=.45\textwidth]{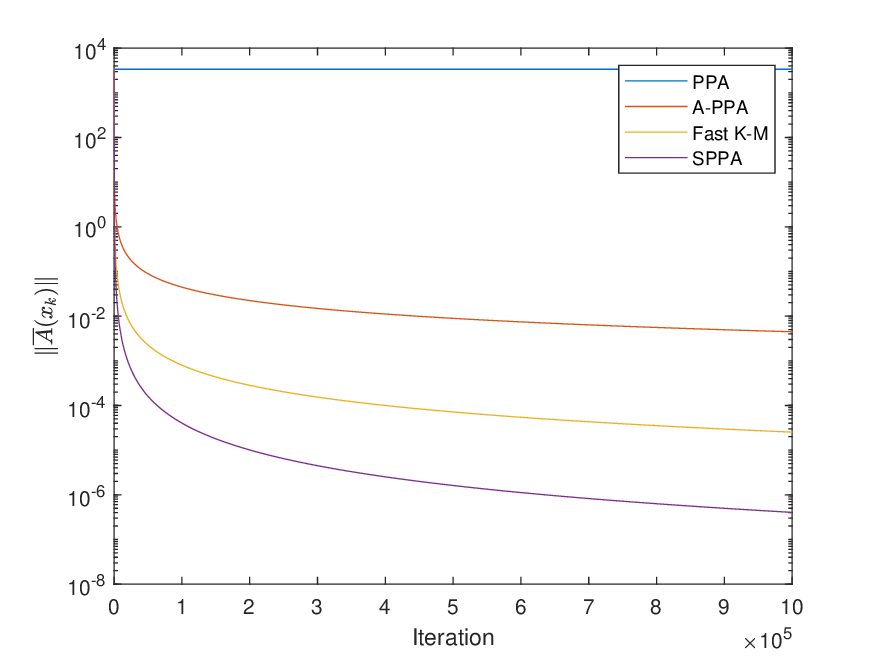}}
	\subfigure{\includegraphics[width=.45\textwidth]{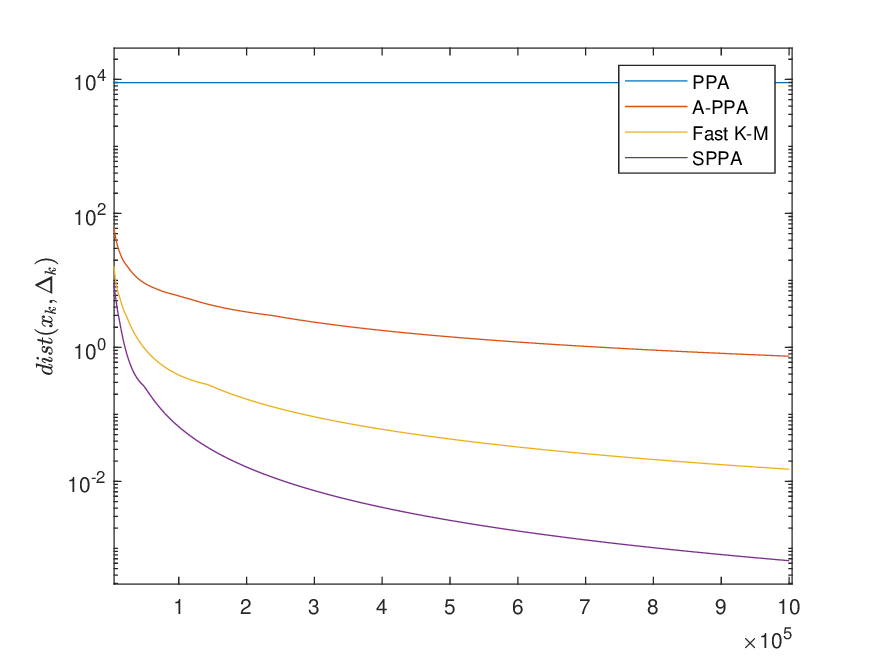}}
	\caption{Comparison between PPA and different accelerated PPAs.}
	\label{fig:parallel}
\end{figure}

In our experiment, \( C_1 = \mathbb{R}_d^{\geqslant 0} \) and \( C_2 = \{ x \in \mathbb{R}^d \mid \langle\vec{1}_d, x\rangle = 1 \} \). It is easy to see that \( C_1 \cap C_2 \) is the unit simplex \( \Delta_d \). Additionally, \( \lambda_k \) is set to 1, \( \omega_1 \) is set to 0.5, and the initial point is randomly generated. We run PPA, accelerated PPA, the Fast K-M algorithm with \( s = 2 \) and \( \alpha = 3 \), and SPPA with \( r = 2 \) and \( C = 1 \) for $10^6$ iterations. Then we plot \( \norm{\overline{A}(x_k)} \) and \( \mathrm{dist}(x_k, \Delta_d) \) with respect to the iteration number, as shown in Figure \ref{fig:parallel}.

As we can see from Figure \ref{fig:parallel}, the PPA converges slowly, whereas the accelerated PPAs converge faster than the PPA. Additionally, our SPPA converges faster than both the accelerated PPA and the Fast K-M algorithm.

Next, we examine whether the rule "the greater \( C \) is, the faster SPPA performs" still holds. Similar to the numerical experiment conducted in Section \ref{sec:num1}, we first compare the numerical performance of SPPAs with \( C = 0.01 \), \( C = 0.25 \), \( C = 0.5 \), \( C = 0.75 \) and \( C = 1 \). We run SPPAs with different values of \( C \) for $10^6$ iterations and plot \( \min_{1 \leqslant k \leqslant 10^6} \norm{\overline{A}(x_k)} \) against \( C \).
\begin{figure}[ht]
	\centering
	\subfigure{\includegraphics[width=.45\textwidth]{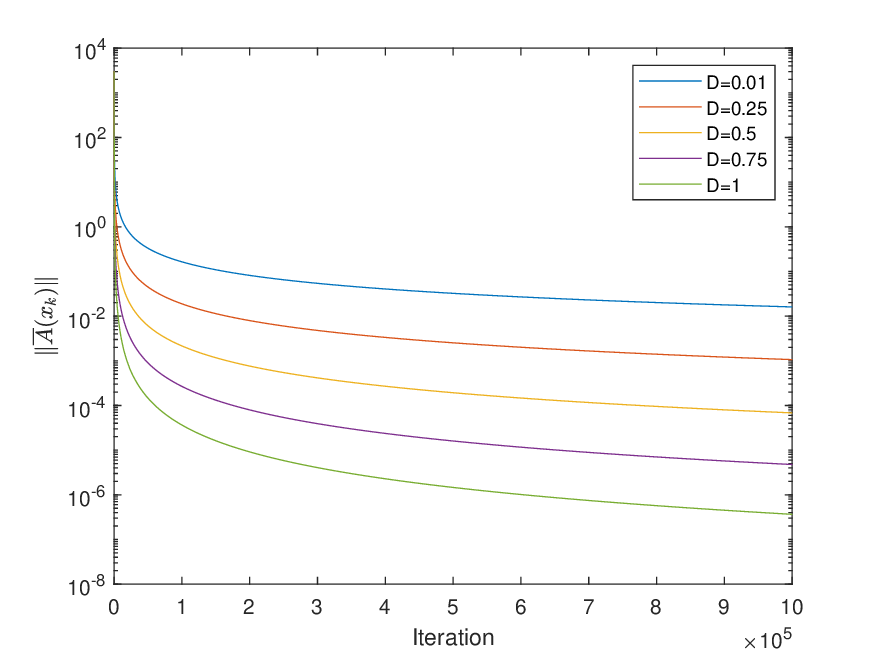}}
	\subfigure{\includegraphics[width=.45\textwidth]{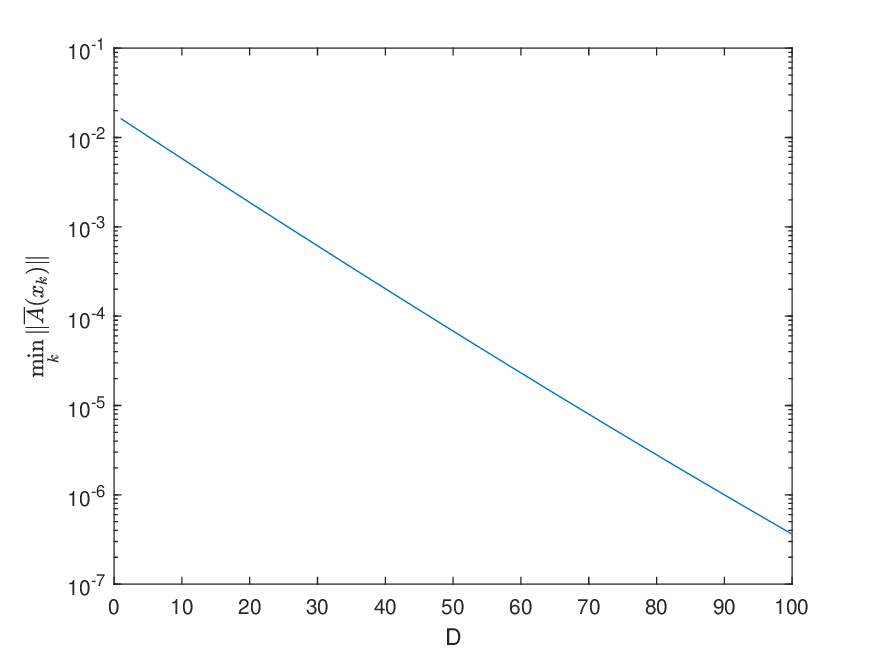}}
	\caption{Comparison between SPPAs with different $C$.}
	\label{fig:parallelD}
\end{figure}

As we can see from Figure \ref{fig:parallelD}, the rule "the greater \( C \) is, the faster SPPA performs" still holds.

Next, we demonstrate that the rule for \( r \) also holds. First, we run SPPAs with \( r = 2 \), \( r = 5 \) and \( r = 10 \), setting \( C = r - 1 \). Then we plot \( \norm{\overline{A}(x_k)} \) against the iteration number. Next, we run SPPAs with different values of \( r \) for \( 10^6 \) iterations and plot $\min_{1 \leqslant k \leqslant 10^6}\norm{\overline{A}(x_k)}$ respect to the iteration number.

\begin{figure}[ht]
	\centering
	\subfigure{\includegraphics[width=.45\textwidth]{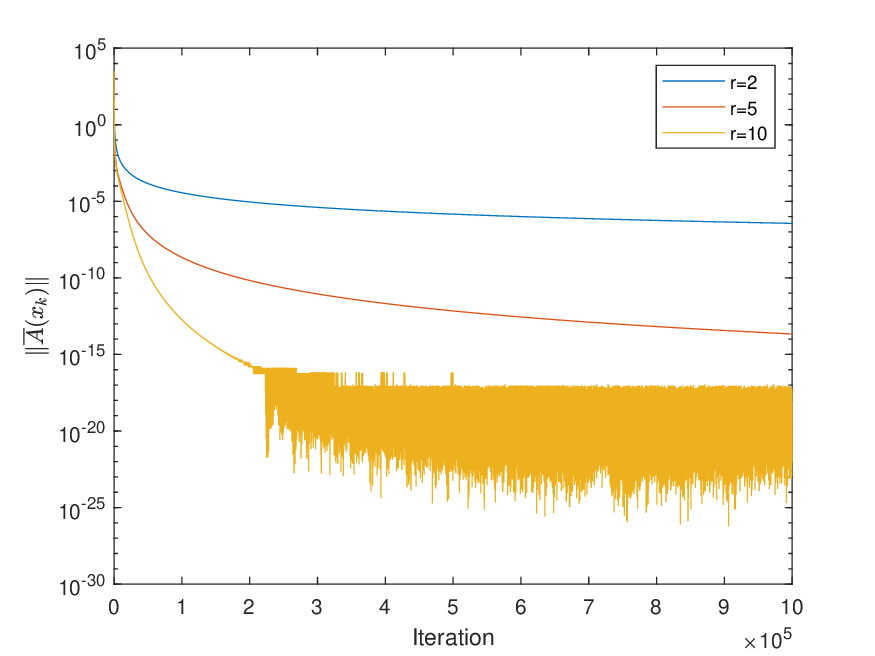}}
	\subfigure{\includegraphics[width=.45\textwidth]{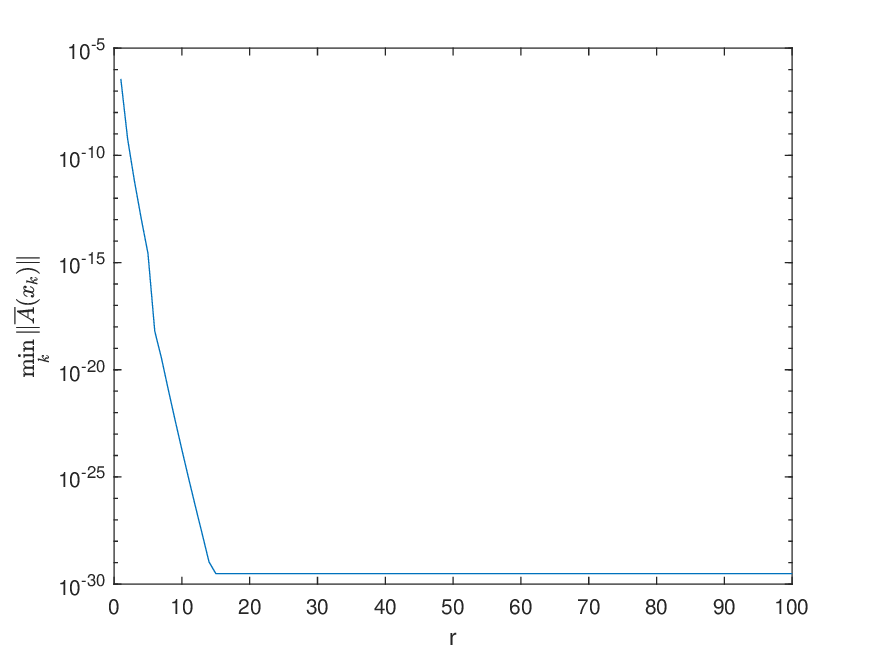}}
	\caption{Comparison between SPPAs with different $r$.}
	\label{fig:parallelR}
\end{figure}

As we can see from Figure \ref{fig:parallelR}, as \( r \) increases, the resulting SPPA converges faster. However, the oscillatory behavior of SPPA with \( r = 10 \) can be observed on the left-hand side of Figure \ref{fig:parallelR}. One possible reason for this phenomenon is that the SPPA has already achieved a solution with very high accuracy so that SPPA suffers from the errors introduced by floating-point arithmetic.

From these numerical experiments, we can conclude the following two rules:
\begin{center}
	The greater \( C \) is, the faster SPPA performs. \\
	Larger \( r \) leads to faster SPPA.
\end{center}

\subsection{Matrix Game Example}

In this section, we consider the min-max problem:
\begin{equation}
	\min_x \max_y f(x) + \inner{Ax, y} - g(y),
	\label{eq:minmax}
\end{equation}
where \( f \) and \( g \) are closed, proper and convex functions defined on \(\mathbb{R}^n\) and \(\mathbb{R}^m\), respectively. The PDHG method, whose iteration rule is given by:
\begin{equation}
	\begin{split}
		x_{k+1} &= \argmin_{x} \left\{ f(x) + \inner{Ax, y_k} + \frac{1}{2\tau} \norm{x - x_k}^2 \right\}, \\
		y_{k+1} &= \argmin_{y} \left\{ g(y) - \inner{A(2x_{k+1} - x_k), y} + \frac{1}{2\sigma} \norm{y - y_k}^2 \right\},
	\end{split}
	\label{eq:PDHG}
\end{equation}
is an effective method to solve \eqref{eq:minmax}. It has been proven that if \( \sigma \tau \norm{A}^2 \leqslant 1 \), \eqref{eq:PDHG} is equivalent to a preconditioned PPA (see \cite{chambolle16,connor20}). The associated preconditioner is:
\[
P = \begin{pmatrix}
	\tau^{-1} \mathbb{I} & -A^\top \\
	-A & \sigma^{-1} \mathbb{I}
\end{pmatrix}.
\]

Due to the fact that \eqref{eq:PDHG} can be regarded as a preconditioned PPA, we test the numerical performance of the symplectic PDHG derived from SPPA. Here, we consider the following matrix game problem:
\begin{equation}
	\min_{x \in \Delta_n} \max_{y \in \Delta_m} y^\top A x,
	\label{eq:matrixgame}
\end{equation}
where \(\Delta_m\) is the \((m-1)\)-dimensional unit simplex and \(A \in \mathbb{R}^{m \times n}\). Problem \eqref{eq:matrixgame} can be transformed into:
\[
\min_x \max_y I_{\Delta_n}(x) + \inner{Ax, y} - I_{\Delta_m}(y),
\]
where \(I_{\Delta_m}\) is the indicator function of \(\Delta_m\). Thus, the PDHG method can be used to solve \eqref{eq:matrixgame}.

In our experiment, \(m = 1000\), \(n = 2000\) and the matrix \(A\) is randomly generated. The initial points for \(y\) and \(x\) are set to \(\vec{1}_{1000}/1000\) and \(\vec{1}_{2000}/2000\), which are the barycenters of \(\Delta_{1000}\) and \(\Delta_{2000}\), respectively. The parameters \(\tau\) and \(\sigma\) are set to \(0.99 / \norm{A}\). We run PDHG, accelerated PDHG (derived from \eqref{eq:acceleratedproximal}), fast PDHG with \(s = 2\) and \(\alpha = 3\) (derived from \eqref{eq:fastkm}), and symplectic PDHG with \(r = 2\) and \(C = 1\) for \(10^5\) iterations. Because PDHG can be viewed as a preconditioned PPA, the term \(\inner{P(\tilde{x}_k - x_k), \tilde{x}_k - x_k}\) plays the same role as \(\norm{\overline{A}(x_k)}^2\) in SPPA. So, we plot \(\inner{P(\tilde{x}_k - x_k), \tilde{x}_k - x_k}\) and the duality gap with respect to the iteration number.

\begin{figure}[ht]
	\centering
	\subfigure{\includegraphics[width=.45\textwidth]{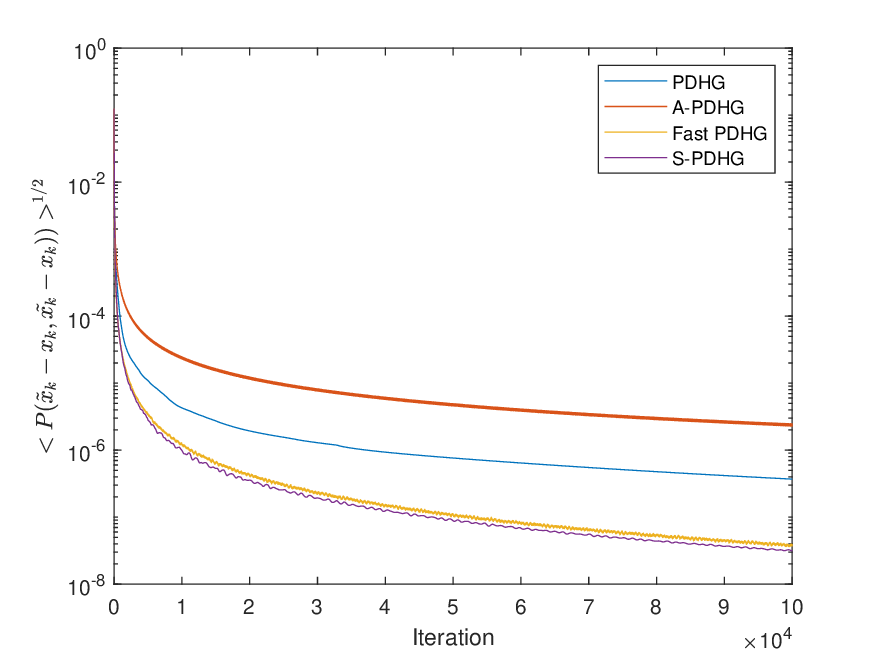}}
	\subfigure{\includegraphics[width=.45\textwidth]{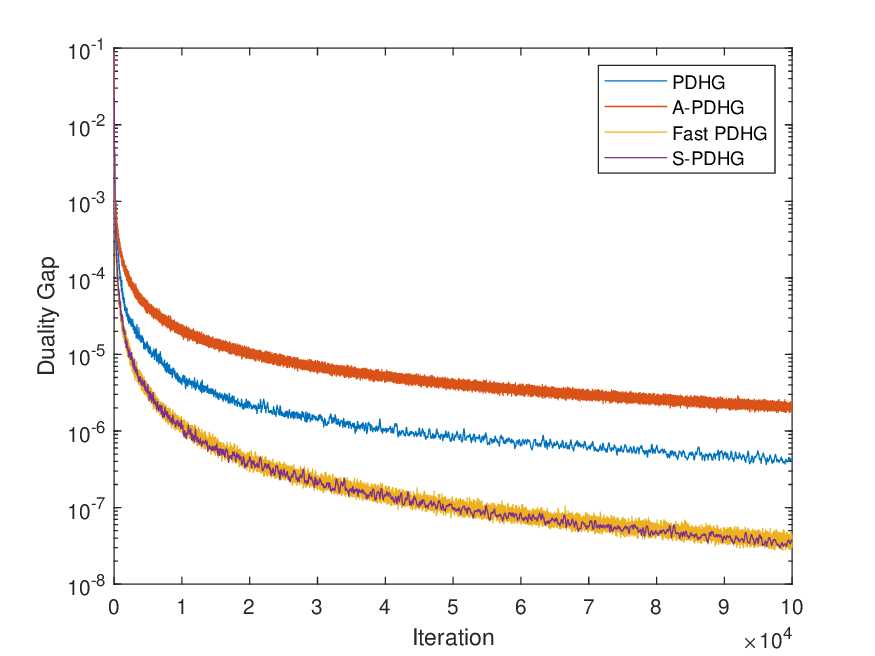}}
	\caption{Comparison between PDHG and several accelerated PDHG methods.}
	\label{fig:PDHG}
\end{figure}

As we can see from Figure \ref{fig:PDHG}, the accelerated PDHG does not perform as well as PDHG, and fast PDHG and symplectic PDHG converge faster than PDHG. Also, the numerical performance of symplectic PDHG is slightly better than fast PDHG on $\inner{P(\tilde{x}_k-x_k),\tilde{x}_k-x_k}$, and the oscillatory behavior of symplectic PDHG is milder than fast PDHG when consider the duality gap.

\section{Conclusion and Discussion}
\label{sec:6}
Throughout this paper, we first propose an ODE system for the maximally monotone operator \( A \). We then employ the Lyapunov function technique and two auxiliary functions to establish that the ODE system exhibits $o(1/t^2)$ convergence rate, and that its solution trajectory converges weakly to a point in $\zero(A)$. Building upon this foundation, we apply the Symplectic Euler Method to discretize the ODE system, leading to the development of our accelerated proximal point algorithm, which we call the Symplectic Proximal Point Algorithm. Furthermore, using the Lyapunov function technique and two auxiliary sequences, we demonstrate that the SPPA also achieves a convergence rate of \( o(1/k^2) \), and the sequences it generates converge weakly to a point in $\zero(A)$.

There are numerous theoretical issues regarding accelerated proximal point algorithms that merit further investigation. First, in \cite{attouch19,chen23}, the authors discussed the convergence rates of the NAG method and FISTA with wider range of parameter. The convergence properties of \eqref{eq:acceleratedode} and Algorithm \ref{al:sppa} when $r<1$ or $C>r-1$ are also topics worthy of exploration. Furthermore, the works in \cite{bao23} and \cite{li24} have shown that the NAG method exhibits exponential convergence under the assumption of strong convexity, even when the strong convexity parameter is unknown. It is valuable to determine whether our SPPA, Fast K-M algorithm or accelerated PPA can achieve (possibly asymptotic) exponential convergence. Additionally, \cite{park22} introduced an accelerated PPA used for \eqref{eq:main} when $A$ is strongly monotone. Therefore, it is pertinent to consider whether we can develop an ODE system for a strongly maximally monotone operator and derive a novel accelerated PPA to solve \eqref{eq:main} under the assumption of strong monotonicity.

\begin{figure}[ht]
	\centering
	\subfigure{\includegraphics[width=.45\textwidth]{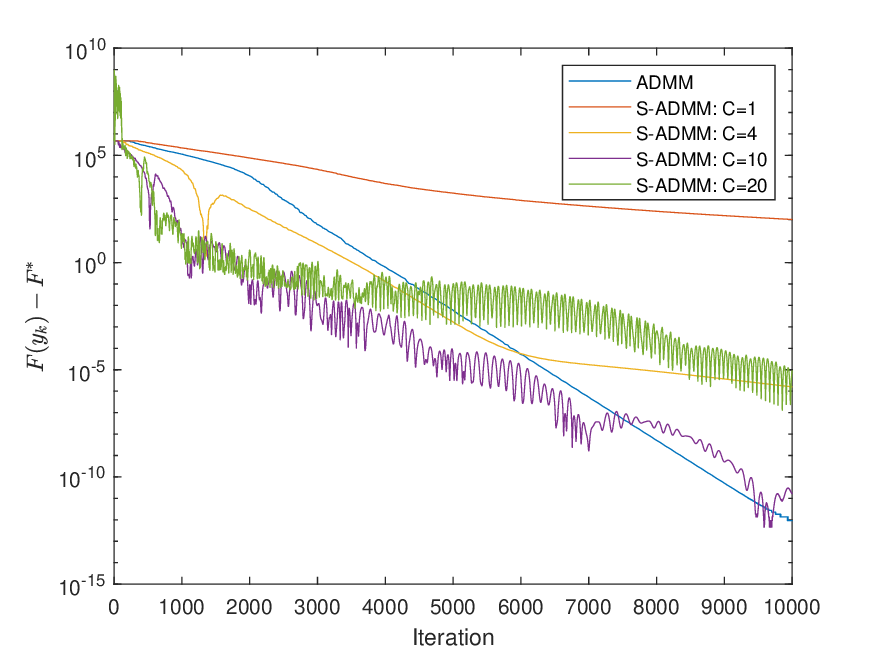}}
	\subfigure{\includegraphics[width=.45\textwidth]{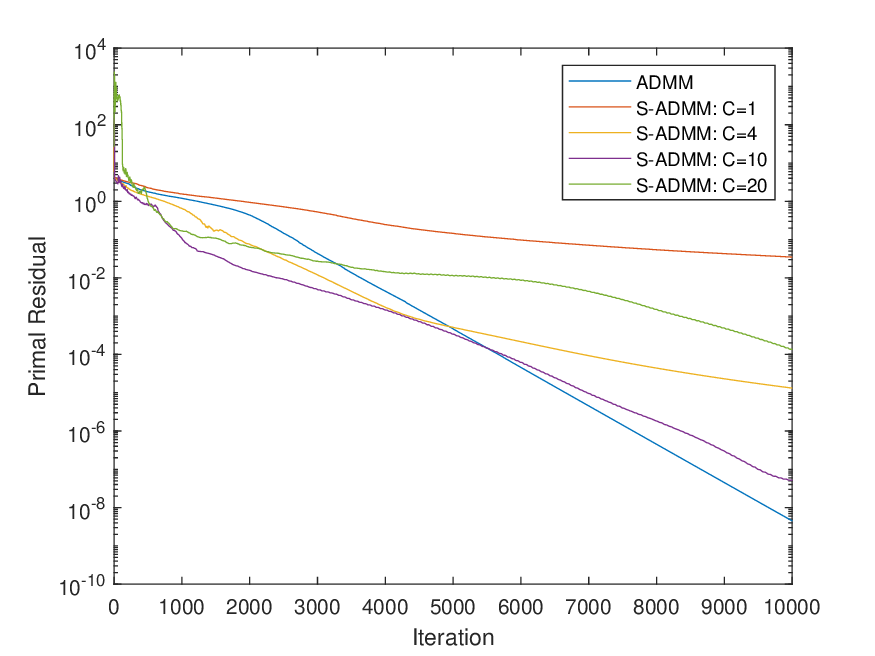}}
	\caption{Evolution of $F(y_k)-F^*$ and primal residual versus iterations on LASSO problem, where $F^*$ is the optimal value of LASSO problem. Here the numerical performance of S-ADMMs with $r=2$ and different $C$ (derived from Algorithm \ref{al:sppa}) is plotted above. It is worthwhile to mention that the convergence rate of ADMM is $O(1/k^2)$ when applied to LASSO problem (see \cite{tian19}). This explains why, regardless of how $r$ is increased in S-ADMM, the corresponding S-ADMM does not converge faster than ADMM. But if we allow $C$ greater than $r-1$, the resulting S-ADMM converges faster in this case. As we can see, S-ADMM with \( C = 10 \) converges faster than ADMM at the beginning. Due to oscillation phenomenon, S-ADMM with \( C = 10 \) eventually converges more slowly than ADMM. However, the parameter \( C \) cannot be arbitrarily large. As evident from the figures, the overall numerical performance of S-ADMM with \( C = 20 \) is not as good as that of S-ADMM with \( C = 10 \).}
	\label{fig:sadmmnumerical3}
\end{figure}

In practice, there are some valuable problems regarding the parameter \( C \). Due to the assumption in Theorem \ref{thm:symplecticpparate}, we have tested SPPA under the condition \( C \leqslant r - 1 \). However, we also observe that for some optimization problems or algorithms, \( C \) can be chosen to be much larger than \( r - 1 \), which can correspond to a faster convergence rate for SPPA. The numerical experiment presented in Figure \ref{fig:sadmmnumerical3} is one such example.

 As we can see from Figure \ref{fig:sadmmnumerical3}, S-ADMM with \( C = 10 \) converges much faster than S-ADMM with \( C = 1 \), but S-ADMM with \( C = 20 \) converges slower than S-ADMM with \( C = 10 \). Thus, what is the optimal value of \( C \) for a specific maximally monotone operator \( A \) that ensures the fastest convergence of SPPA? Moreover, determining the optimal value of \( C \) for a wide variety of maximally monotone operators can be very time-consuming and resource-intensive. Developing an algorithm that yields a sufficiently large value of \( C \) to maintain the accelerated behavior of Algorithm \ref{al:sppa} is crucial for practical applications.

\bibliographystyle{amsplain}
\bibliography{reference}

\end{document}